\documentclass[a4paper,12pt,leqno]{amsart}

\usepackage{amsmath}
 \usepackage{amssymb}
\usepackage{amsthm}
\usepackage{amsopn}
\usepackage{amscd}
\usepackage{xcolor}
\usepackage{tikz}
\usepackage[utf8]{inputenc}
\usepackage[T1]{fontenc}
\usepackage[colorlinks=true,linkcolor=blue,citecolor=blue,urlcolor=blue]{hyperref}
\newcommand{\frg}{\mathfrak{g}}

\newcommand{\frmm}{\mathfrak{m}}
\newcommand{\frn}{\mathfrak{n}}

\newcommand{\bbC}{\mathbb{C}}

\newcommand{\bbR}{\mathbb{R}}

\newcommand{\bbZ}{\mathbb{Z}}

\newcommand{\caI}{\mathcal{I}}
\newcommand{\caJ}{\mathcal{J}}
\newcommand{\caK}{\mathcal{K}}

\newcommand{\caM}{\mathcal{M}}

\newcommand{\caO}{\mathcal{O}}

\newcommand{\caR}{\mathcal{R}}

\theoremstyle{plain}
\newtheorem{thm}{Theorem}[section]
\newtheorem{lemma}[thm]{Lemma}
\newtheorem{cor}[thm]{Corollary}
\newtheorem{prop}[thm]{Proposition}

\theoremstyle{definition}
 \newtheorem{defi}[thm]{Definition}

\theoremstyle{definition}

\newtheorem{rmk}[thm]{Remark}
\newtheorem{example}[thm]{Example}

\newtheorem{thmA}{Theorem}

\newtheorem{thmLM}{Theorem}

\newcommand{\Id}{\mathrm{Id}}

\newcommand{\GL}{\mathbf{GL}}

\begin{document}

\numberwithin{equation}{section}

\title{Resolution of ladder representations}
\author{David Renard}
\address{
Centre de Mathématiques \\
Laurent Schwartz
\\ Ecole Polytechnique }
\date{\today}

\begin{abstract}
Let $\pi$ be a ladder representation of $\mathbf{GL}_m(F)$, $F$ a non-Archimedean local field. We
construct an exact sequence of costandard modules, indexed by a lower Bruhat interval of a symmetric
group and graded by the length, whose Euler characteristic is the determinantal formula of Tadi\'c and
Lapid--M\'inguez. This is an analogue of the dual Bernstein--Gelfand--Gelfand resolution. The proof is
combinatorial, through Jacquet modules and the Bruhat order; we also interpret the coresolution as a
Cousin complex on Zelevinsky's space, which extends it to the balanced representations of
Lapid--M\'inguez.
\end{abstract}

\maketitle

%\frontmatter
%\tableofcontents

%\mainmatter

\section*{ Introduction}
Let $F$ be a non-Archimedean locally compact  field and, for every integer $m\geq0$, set $G_m=\GL(m,F)$
(with the convention that $G_0$  is the trivial group). Denote by $\caR_m$ the Grothendieck group
of the category of
 smooth representations of $G_m$  of  finite length, and let $\caR=\oplus_{m\geq 0}\caR_m$ with the
 product structure
 defined by normalized  parabolic induction. The commutative ring $\caR$  (with the one-dimensional
 representation
 of $G_0$ as the identity element) was introduced by Zelevinsky, who showed that it is freely generated by the
 essentially square-integrable representations   (of any $G_m$) \cite{Zel80}. Hence, the monomials in these
generators (the so-called standard modules) form a basis for   $\caR$. Another natural basis for $\caR$ is
given by the irreducible representations. The change of basis matrix is unitriangular (in an
appropriate sense)
 and the coefficients can be expressed in terms of the Kazhdan-Lusztig polynomials.
This fact was conjectured by Zelevinsky \cite{Zel81} (and in a more precise form in \cite{Zel85})
and proved in \cite{CG97}.
 This solves in principle the problem of decomposing a standard module into irreducible constituents, 
 and, by inversion, of expressing an irreducible representation in terms of standard modules.
 However, Kazhdan-Lusztig 
 polynomials are computed by a recursive algorithm, and they do not give closed formulas in
 general. Nevertheless, 
 in the case of a Speh representation, Tadi\'c obtained a remarkable closed  formula expressing its character 
 as a linear combination of characters of standard modules with coefficients $\pm1$  \cite{Tad95}.  
 Tadi\' c ingenious argument transfers the problem to a question about complex
 groups. Later on, the argument was  simplified by Chenevier and Renard \cite{CR08} by the use of
 the Desnanot-Jacobi
identity  of determinants (also known as Dodgson's rule of determinants or Lewis Carroll identity).
Both proofs rely on the
determination of the composition series at the edge of a complementary series, and thus they
heavily rely on unitarity.
In \cite{LM13}, Lapid and Minguez give a different proof of Tadic's formula
which also describes the maximal proper submodule of the standard module, that is, the kernel of
its projection onto the Langlands quotient (Theorem \ref{LM} (i) below), and at the same time
extends to a wide class of representations called ladder  representations.
Except for the Speh representations, ladder representations are not unitary. Instead of
unitarity, they use Jacquet module technique. They actually need to  study not only the
semisimplification of the Jacquet
module, but the finer structure of its submodules.

To state the determinantal formula of Tadic, let us introduce some more notation. Denote by  $\nu$
the character
$\vert\det \vert$  on any $G_n$. (The $n$ will be implicit and hopefully clear from the context.) For any
smooth representation $\pi$ of $G_n$  and any $a\in \bbR$ denote by $\nu^a  \pi$ the representation
obtained from $\pi$
by twisting it by the character $\nu^a$. If $\pi_1, \ldots ,  \pi_r$ are  smooth representations of
 $G_{n_1}, \ldots ,G_{n_r}$ respectively, we will denote as usual by
\[\pi_1\times \ldots \times \pi_r\]
the representation of $ G_{n_1+\cdots +n_r}$ parabolically induced (normalized induction) from the
representation
$\pi_1\otimes \ldots \otimes \pi_r$ of the standard parabolic subgroup of $G_{n_1+\cdots +n_r}$ of
type $(n_1,\ldots, n_r)$ (with Levi subgroup   $G_{n_1}\times \cdots \times G_{n_r}$).
Throughout we  fix a positive integer $d$ and an irreducible cuspidal (not necessarily unitary) representation
$\rho$ of $G_d$.  By a segment $[a;b]$  we mean a set of $\{a,a+1,\ldots , b\}$ where $a\leq b$ are integers.
For any segment   $[a;b]$   we denote by $\Delta([a;b])$  the unique irreducible quotient of
 $\nu^a\rho\times \cdots \times \nu^b\rho$.
 It is an essentially square-integrable representation of $G_{(b-a+1)d}$. The map $[a;b]\mapsto
 \Delta([a;b])$ is a bijection
 between the segments and the essentially square-integrable representations whose cuspidal support
 is contained
in $\{ \nu^i \rho: i\in \bbZ\}$
\cite{Zel80}, \S[9.3]. We also write $\Delta([a; a-1]) = 1$ (the one-dimensional representation of $G_0$)
and for any $a\in \bbZ$, $\Delta([a; b]) = 0$  if $b < a-1$.
Let   $[a_1;b_1],\ldots ,[a_t;b_t]$  be segments and let $m_i=b_i-a_i+1$ and $\Delta_i=\Delta([a_i;b_i])$.
 Assume that $\Delta_i$  does not precede $\Delta_j$
for any $i<j$, i.e., we do not have  $a_i<a_j\leq b_i+1<b_j+1$ for $i<j$. Then the
representation  $\Delta_1\times \cdots \times \Delta_t$ admits a unique irreducible quotient, the
Langlands quotient,
which we denote by $L(\Delta_1,\ldots ,\Delta_t) $ (see for example \cite{Rod82}, Th{\' e}or{\`e}me 3).
Assume now that $a_1>\ldots  >a_t$ and $b_1>\cdots >b_t$. In this case we say that
$L(\Delta_1,\ldots ,\Delta_t) $ is a
 ladder representation.
 A particularly important subclass is the Speh representations where for every $i=1,\ldots,t-1$,
 $a_{i+1}=a_i-1$ and $b_{i+1}=b_i-1$.
These representations comprise the building blocks for the unitary dual of $G_n$ (\cite{Tad86}).
Let $S_t$  be the symmetric group on $\{1,\ldots, t\}$. For any  $w\in S_t$ let
\[\caI_w = \Delta([a_{w(1)}; b_1])\times \cdots \times \Delta([a_{w(t)}; b_t]).\]
(This is $0$ unless $a_{w(i)}\leq b_i+1$ for all $i$.) In particular, for $k=1, \ldots ,t-1$ let $\caK_k =
 \caI_{s_k}\subset \caI_{\Id}$
where $s_k$  is the transposition $(k,k + 1)$.

\begin{thmLM}[Lapid, M\'inguez]\label{LM} Assume that $L(\Delta_1,\ldots ,\Delta_t) $ is a ladder
representation
and let $\pi=L(\Delta_1,\ldots ,\Delta_t) $. Then

(i) The maximal proper submodule of $\Delta_1\times \cdots \times \Delta_t$ is $\caK=\sum_{k=1}^{t-1}\caK_k$.

(ii) $\pi$ satisfies Tadic's determinantal formula (in $\caR$)
\begin{equation} \pi=\sum_{w\in S_t} (-1)^{l(w)} \caI_w
\end{equation}
\end{thmLM}

Tadic's determinantal formula is reminiscent of Weyl character formula for finite dimensional
representations of a complex,  finite dimensional,   semisimple Lie algebra $\frg$. For any integral weight
$\lambda$, denote by  $M(\lambda)$ the Verma module with highest weight $\lambda$ and $L(\lambda)$
its unique irreducible subquotient.
If $\lambda$ is dominant, then $L(\lambda)$ is finite dimensional, and in that case
 \begin{equation*} L(\lambda)=\sum_{w\in W} (-1)^{l(w)} M(w\cdot \lambda)
\end{equation*}
in the Grothendieck group of the category $\caO$.
The Weyl character formula for $L(\lambda)$ is the trace in the Grothendieck group of a much deeper fact,
namely the existence  of the Bernstein-Gelfand-Gelfand resolution in category $\caO$. This is a resolution of
the finite dimensional irreducible module $L(\lambda)$ by direct sums of  Verma modules $M(w\cdot \lambda)$.
 More precisely, the resolution is of the form
\[ 0\rightarrow  M(w_0\cdot \lambda) \rightarrow \bigoplus_{w\in W\vert l(w)=l(w_0)-1} M(w\cdot \lambda)
 \rightarrow\cdots \qquad \qquad\qquad \qquad \]
 \[ \qquad \qquad\qquad \qquad\cdots  \rightarrow \bigoplus_{w\in W\vert l(w)=1}   M(w\cdot
 \lambda) \rightarrow
M(\lambda)\rightarrow L(\lambda)\rightarrow 0
\]
where $w_0$ is the longest element in the Weyl group $W$ of $\frg$.
The Euler-Poincaré principle recovers Weyl character formula from the  Bernstein-Gelfand-Gelfand resolution.

Theorem \ref{LM} strongly suggests that a similar resolution of a ladder representation
$\pi=L(\Delta_1,\ldots ,\Delta_t) $
by direct sums of standard representations $\caI_w$  should exist. The
Bernstein-Gelfand-Gelfand resolution has a geometric construction (Kempf, Mili\v{c}i\'c): the Cousin
complex of the structure sheaf of the flag variety, with respect to the stratification by the Bruhat
cells, gives after taking global sections the exact sequence
\[0\rightarrow L(\lambda)\rightarrow M(\lambda)^\vee\rightarrow\bigoplus_{l(w)=1}M(w\cdot\lambda)^\vee
\rightarrow\cdots\rightarrow M(w_0\cdot\lambda)^\vee\rightarrow0,\]
where $M(\mu)^\vee$ is the dual Verma module; it is the dual (contragredient) of the BGG resolution
(recall that $L(\lambda)$ is self-dual, \cite{Hum08}, \S3.2 and Theorem 3.3(a)),
and the term of degree $k$ comes from the cells of codimension $k$ (\S\ref{ssec:geomint}). In our setting,
 the analogues of the dual Verma modules are the costandard modules
\[\caJ_w = \Delta([a_{w(t)}; b_t])\times \cdots \times \Delta([a_{w(1)}; b_1]), \qquad w\in S_t,\]
in which the segments appear in the order opposite to that of $\caI_w$. The representation $\pi$
 is the unique irreducible subrepresentation of $\caJ_{\Id}$, and $\caJ_w$ and $\caI_w$ have
 the same image in $\caR$. In this paper, we prove the existence of such a dual
 Bernstein-Gelfand-Gelfand resolution of a ladder representation.

Such a resolution can also be obtained, for modules over the graded Hecke algebra of type $A$,
by applying the Arakawa-Suzuki functor to the BGG resolution of finite-dimensional
$\mathfrak{gl}_n$-modules (\cite{BC15}, Remark 4.3.2). Our proof, by contrast, takes place entirely 
in the category of smooth representations of the $p$-adic groups $G_n$. 
Apart from Theorem \ref{LM} (ii) and the Langlands classification \cite{Rod82}, the proof only uses
the results of Zelevinsky 
and Bernstein-Zelevinsky \cite{BZ77, Zel80}, the computation of Jacquet modules of ladder representations 
by Kret and Lapid \cite{KL12}, and the combinatorics of the Bruhat order. In particular, it does
not use Kazhdan-Lusztig polynomials.
 The relation with the Kazhdan-Lusztig polynomials is explained in \S\ref{ssec:KL}.

Let $S_t^0=\{w\in S_t : a_{w(i)}\leq b_i+1,\ i=1,\ldots,t\}$, so that $\caI_w\neq 0$ (equivalently
$\caJ_w\neq 0$)
if and only if $w\in S_t^0$. By \cite{LM13}, Lemma 15, $S_t^0$ is the Bruhat interval $\{w\in S_t :
w\leq w'_0\}$
 for a unique $w'_0\in S_t$ (see Section~\ref{sec:bruhat}).

\begin{thmA}\label{BGG}
 There exists an exact sequence of the form
\[ 0\longrightarrow \pi \stackrel{\delta_0}{\longrightarrow} \caJ_{\Id}
\stackrel{\delta_1}{\longrightarrow} \bigoplus_{w\in S_t^0 \vert l(w)=1} \caJ_w
\stackrel{\delta_2}{ \longrightarrow} \cdots  \qquad \qquad \]
\[ \qquad \qquad \cdots   \longrightarrow \bigoplus_{w\in S_t^0\vert l(w)=l(w'_0)-1}  \caJ_w
  \stackrel{\delta_{l(w'_0)}}{\longrightarrow} \caJ_{w'_0}\longrightarrow 0.
\]
\end{thmA}

The contragredient of a ladder representation is again a ladder representation:
 $\pi^\vee=L(\Delta'_1,\ldots,\Delta'_t)$ with $\Delta'_i=\Delta^\vee_{t+1-i}$, i.e. $a'_i=-b_{t+1-i}$
 and $b'_i=-a_{t+1-i}$ (with $\rho$ replaced by $\rho^\vee$). Since taking contragredients commutes
  with normalized parabolic induction, one checks that $\caJ_w^\vee\simeq \caI'_{w_0w^{-1}w_0}$,
  where $\caI'$ denotes the standard modules attached to $\pi^\vee$ and $w_0$ is the longest
  element of $S_t$; the map $w\mapsto w_0w^{-1}w_0$ preserves the length and the Bruhat order.
  Hence Theorem \ref{BGG} is equivalent to the existence of a resolution
\[0\rightarrow \caI'_{w'_0}\rightarrow \cdots\rightarrow \bigoplus_{l(w)=1}\caI'_w\rightarrow
\caI'_{\Id}\rightarrow \pi^\vee\rightarrow 0\]
of the ladder representation $\pi^\vee$ by standard modules. In particular, the first step of the
exact sequence of Theorem \ref{BGG} is the contragredient of Theorem \ref{LM} (i) applied to $\pi^\vee$.

Both Theorem \ref{LM} (ii) and Theorem \ref{BGG} extend to \emph{generalized ladders}, in which
empty segments 
are allowed (Lemma \ref{lem:empty}); these occur naturally as minors of ladders, and are needed in
the induction.
 As a by-product of the proof, we also obtain a determinantal formula of Tadi\'c type for the representations 
 $L(\caM_{c_k})$ attached to the minors, which are not ladder representations in general (Remark
 \ref{rmk:tadic_ck}). Finally, Theorem \ref{BGG} extends to the \emph{balanced} representations of
 Lapid and M\'inguez \cite{LM18}, a class which contains the ladder representations and is
 characterized, among the representations attached to regular multisegments, by the irreducibility
 of $\pi\times\pi$; we sketch this extension, which follows from
 the geometric interpretation of the coresolution as a Cousin complex, in \S\ref{ssec:balanced}.

Atobe \cite{Ato24} has introduced an analogue of ladder representations for the split groups
$\mathbf{Sp}(2n,F)$ and $\mathbf{SO}(2n+1,F)$, which contains the strongly positive discrete series and
the irreducible representations with irreducible $A$-parameters, and has proved a formula expressing
their classes in the Grothendieck group as alternating sums of classes of standard modules.
 In \cite{Ren26b}, we prove that this formula is, after cancellation, the Euler characteristic of a coresolution
by costandard modules analogous to Theorem \ref{BGG}, graded by the dimensions of the orbits of the
Langlands parameters in the graded Lie algebra attached to the infinitesimal character. In particular all the coefficients
of Atobe's formula are $\pm1$, as Atobe expected. Theorem \ref{BGG} is graded in the same way: for
$\GL_n$, the graded Lie algebra attached to the infinitesimal character is Zelevinsky's space of graded
nilpotent endomorphisms, and the Bruhat length $l(w)$ of the permutation indexing a term $\caJ_w$ is the
dimension of the orbit of its Langlands parameter, relative to the orbit of the parameter of $\pi$
(Lemma \ref{lem:dimorb}).

\subsection*{Resolutions of irreducible representations by standard modules}
Theorem \ref{BGG} belongs to a family of results in which an irreducible representation is resolved
by standard (or costandard) modules, one for each stratum of a geometric stratification, placed in
the degree given by the codimension of the stratum:
\begin{itemize}
\item the Bernstein-Gelfand-Gelfand resolution of a finite-dimensional representation by Verma
modules \cite{BGG75}, whose dual is the Cousin resolution of the structure sheaf of the flag variety
along the Bruhat cells (Kempf \cite{Kem78}, Mili\v{c}i\'c \cite{Mil});
\item Johnson's resolutions of the cohomologically induced modules $A_\mathfrak{q}(\lambda)$ by
standard Harish-Chandra modules \cite{Joh84}, graded by the dimensions of $K$-orbits, and the
resolutions of $\theta$-twisted Speh representations of \cite{AMR18}, \S10.3, graded by the
$\theta$-length;
\item for $p$-adic $\GL_n$, the image of the dual BGG resolution under the Arakawa-Suzuki functor
\cite{BC15}, and Theorem \ref{BGG};
\item the ladder representations of Atobe for the classical groups \cite{Ren26b}, graded by the dimensions of the orbits of the Langlands parameters.
\end{itemize}
In all these cases, the mechanism is that of a Cousin resolution. There is a stratified variety,
and the irreducible representation corresponds to a sheaf (on the side exchanged with the original
one by the Fourier transform, or by Koszul duality) whose costalks are concentrated in a single
degree along each stratum, the degree given by its codimension, as for the constant sheaf of a
(rationally) smooth variety. The Cousin complex of this sheaf with respect to the stratification then
gives the coresolution, its terms being the costandard objects ($*$-extensions of the local systems
on the strata). In a highest weight category such as $\mathcal O$, this condition says that the
$\mathrm{Ext}$ groups of the irreducible module with the standard objects are one-dimensional and
concentrated in that degree. For Harish-Chandra modules (Johnson's case, $K$-orbits on the flag
variety), the stabilizers need not be connected, an orbit may carry several local systems, and the
$\mathrm{Ext}$ groups between objects attached to the same orbit involve its equivariant cohomology;
the condition is then expressed on the costalks. 
For ladder representations of $\GL_n$, the stratified space is Zelevinsky's space of graded
nilpotent endomorphisms, stratified by the orbits of the multisegments \cite{Zel81, Zel85}. Along the
strata of a ladder, the Kazhdan-Lusztig polynomials of these orbit closures are those of $S_t$
(Proposition \ref{fact:KL}), and the condition on the costalks becomes $P_{x,w_0}=1$ for all $x$: it
reflects the smoothness of the full flag variety of $\GL_t$ (\S\ref{ssec:geomint}). The present paper and \cite{Ren26b} are part of a
project which aims at understanding the extent of this phenomenon: which irreducible representations
of real or $p$-adic reductive  groups admit such resolutions, and how the grading is determined by the geometry of the
Langlands (or Arthur)  parameters.

\subsection*{On the use of artificial intelligence}
This paper was written with the help of the AI assistant Claude (Anthropic). The starting point is a
draft by the author (2020), which contained the problem, the statements of Theorems \ref{LM} and
\ref{BGG}, the introduction, the preliminaries of \S1, the idea of a proof modelled on the Cousin
resolution, and the beginning of \S2. The assistant corrected several errors of this draft, wrote the
proofs of \S\S2--3 (inclusions of standard modules, minors, key lemma, multiplicity one, construction
of the complex and end of the proof), the comparison with Kazhdan-Lusztig polynomials and the
geometric interpretation of \S\S\ref{ssec:KL}--\ref{ssec:geomint}, the running example and its figures,
and the bibliography, and checked numerically several intermediate statements. The author directed
the work (choice of the strategy, of the references, in particular Johnson's and the Cousin point of
view, and of the presentation) and revised the text. The author takes full responsibility for its
content.

\section{Notation and preliminaries}
\subsection{Conventions}
As in the introduction, $F$ is a non-Archimedean locally compact field, $G_m=\GL(m,F)$, and
$\tau_1\times\cdots\times\tau_r$ is the representation normalized parabolically induced from the
standard (block upper triangular) parabolic subgroup of type $(n_1,\ldots,n_r)$. For a representation
$\tau$ of finite length, we write $[\tau]$ for its image in the Grothendieck ring $\caR$. All the
representations considered have their cuspidal support contained in $\{\nu^i\rho:i\in\bbZ\}$, for the
cuspidal representation $\rho$ of $G_d$ fixed in the introduction. The only Jacquet functor we use is
the minimal one, $r_{\min}$, through the consequence of the geometric lemma of Bernstein and Zelevinsky
(\cite{BZ77}, 2.12) recalled in \S\ref{ssec:words}.

\subsection{Langlands classification}
We use the notation of the introduction for the segments $[a;b]$, the representations
$\Delta([a;b])$ (with the conventions for $b=a-1$ and $b<a-1$), the Langlands quotients
$L(\Delta_1,\ldots,\Delta_t)$, the ladders and the Speh representations. Every irreducible representation with cuspidal support in $\{\nu^i\rho\}$ is of the form $L(\Delta_1,\ldots,\Delta_t)$, for
segments which are uniquely determined up to permutation (\cite{Rod82}, Th\'eor\`eme 3); by reordering
the segments if necessary, $L(\Delta_1,\ldots,\Delta_t)$ is defined for any $t$-tuple of segments. The
segments of a ladder are always numbered so that $a_1>\cdots>a_t$ and $b_1>\cdots>b_t$.

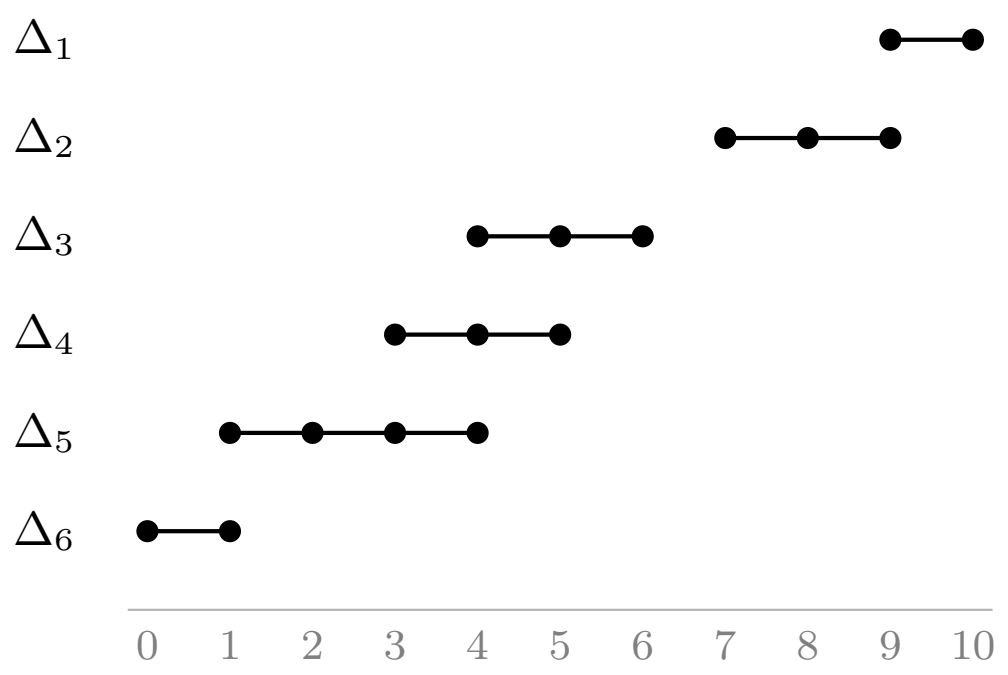
\begin{figure}[h]
\centering
\begin{tikzpicture}[scale=1.3]
\draw[black,thick] (3.78,-0.50)--(4.20,-0.50);
\fill[black] (3.78,-0.50) circle (1.6pt);
\fill[black] (4.20,-0.50) circle (1.6pt);
\node[left,font=\scriptsize] at (-0.25,-0.50) {$\Delta_{1}$};
\draw[black,thick] (2.94,-1.00)--(3.78,-1.00);
\fill[black] (2.94,-1.00) circle (1.6pt);
\fill[black] (3.36,-1.00) circle (1.6pt);
\fill[black] (3.78,-1.00) circle (1.6pt);
\node[left,font=\scriptsize] at (-0.25,-1.00) {$\Delta_{2}$};
\draw[black,thick] (1.68,-1.50)--(2.52,-1.50);
\fill[black] (1.68,-1.50) circle (1.6pt);
\fill[black] (2.10,-1.50) circle (1.6pt);
\fill[black] (2.52,-1.50) circle (1.6pt);
\node[left,font=\scriptsize] at (-0.25,-1.50) {$\Delta_{3}$};
\draw[black,thick] (1.26,-2.00)--(2.10,-2.00);
\fill[black] (1.26,-2.00) circle (1.6pt);
\fill[black] (1.68,-2.00) circle (1.6pt);
\fill[black] (2.10,-2.00) circle (1.6pt);
\node[left,font=\scriptsize] at (-0.25,-2.00) {$\Delta_{4}$};
\draw[black,thick] (0.42,-2.50)--(1.68,-2.50);
\fill[black] (0.42,-2.50) circle (1.6pt);
\fill[black] (0.84,-2.50) circle (1.6pt);
\fill[black] (1.26,-2.50) circle (1.6pt);
\fill[black] (1.68,-2.50) circle (1.6pt);
\node[left,font=\scriptsize] at (-0.25,-2.50) {$\Delta_{5}$};
\draw[black,thick] (0.00,-3.00)--(0.42,-3.00);
\fill[black] (0.00,-3.00) circle (1.6pt);
\fill[black] (0.42,-3.00) circle (1.6pt);
\node[left,font=\scriptsize] at (-0.25,-3.00) {$\Delta_{6}$};
\draw[gray!60] (-0.10,-3.40)--(4.30,-3.40);
\node[font=\tiny,gray] at (0.00,-3.58) {0};
\node[font=\tiny,gray] at (0.42,-3.58) {1};
\node[font=\tiny,gray] at (0.84,-3.58) {2};
\node[font=\tiny,gray] at (1.26,-3.58) {3};
\node[font=\tiny,gray] at (1.68,-3.58) {4};
\node[font=\tiny,gray] at (2.10,-3.58) {5};
\node[font=\tiny,gray] at (2.52,-3.58) {6};
\node[font=\tiny,gray] at (2.94,-3.58) {7};
\node[font=\tiny,gray] at (3.36,-3.58) {8};
\node[font=\tiny,gray] at (3.78,-3.58) {9};
\node[font=\tiny,gray] at (4.20,-3.58) {10};
\end{tikzpicture}
\caption{The ladder of Example \ref{ex:run}.}\label{fig:ladder}
\end{figure}

\begin{example}[The running example]\label{ex:run} Throughout the paper we follow the ladder with
$t=6$ segments
\begin{gather*}[a_1;b_1]=[9;10],\ [a_2;b_2]=[7;9],\ [a_3;b_3]=[4;6],\\
[a_4;b_4]=[3;5],\ [a_5;b_5]=[1;4],\ [a_6;b_6]=[0;1],\end{gather*}
so that $\pi=L(\Delta_1,\ldots,\Delta_6)$ is a representation of $G_{17d}$. In Figure
\ref{fig:ladder}, each segment $[a_i;b_i]$ is drawn as a row of dots at the abscissae
$a_i,\ldots,b_i$ (the abscissa $c$ standing for $\nu^c\rho$); the beginnings and the ends both
decrease from top to bottom. It is not a Speh representation, and consecutive segments are linked in
various ways: $\Delta_1,\Delta_2$ and $\Delta_5,\Delta_6$ overlap in one point, $\Delta_2,\Delta_3$
are juxtaposed, and $\Delta_3,\Delta_4$ and $\Delta_4,\Delta_5$ overlap in two points.
\end{example}

\subsection{Linked segments}\label{ssec:Z}
Recall that two segments are \emph{linked} if neither of them contains the other and their union is
a segment (\cite{Zel80}, 4.1); they are \emph{juxtaposed} if moreover their intersection is
empty, i.e. one of them begins just after the end of the other. An \emph{elementary operation} on a multisegment replaces two linked
segments by their union and their intersection, the latter being omitted if it is empty. For
multisegments $\frn$, $\frmm$, we write $\frn\leq\frmm$ if $\frn$ is obtained from $\frmm$ by a
sequence of elementary operations (\cite{Zel80}, 7.1). We use two facts.
\begin{itemize}
\item[(Z)] (\cite{Zel80}, Theorem 9.7 (a)) A product $\Delta_1\times\cdots\times\Delta_r$ of
essentially square-integrable representations is irreducible if and only if no two of the segments
$\Delta_1,\ldots,\Delta_r$ are linked. In that case, the isomorphism class of the product does not
depend on the order of the factors, since its image in $\caR$ does not.
\item[(L)] (Langlands classification) Let $\frn$ be a multisegment, and let
$\Delta''_1,\ldots,\Delta''_r$ be its segments, ordered so that their ends are non-decreasing. Then
$L(\frn)$ is the unique irreducible subrepresentation of
$\Delta''_1\times\cdots\times\Delta''_r$. This is the well-known subrepresentation form of the
Langlands classification; it follows from the quotient form (\cite{Rod82}, Th\'eor\`eme 3) by taking
contragredients, since $(\Delta_1\times\cdots\times\Delta_r)^\vee\simeq\Delta_1^\vee\times\cdots\times
\Delta_r^\vee$, $\Delta([a;b])^\vee=\Delta([-b;-a])$ (with $\rho$ replaced by $\rho^\vee$), and
$L(\frn)^\vee\simeq L(\frn^\vee)$.
\end{itemize}

\section{Segments and Bruhat order}\label{sec:bruhat}

 Let $ [a_1; b_1], \ldots , [a_t; b_t ]$  be segments forming a ladder.

Let
\[ S_t^0=\{w \in S_t \, \vert \, a_{w(i)}\leq b_i+1, \, i=1,\ldots ,t\}.\]
Note that if $w \notin  S_t^0$, $\caI_w=0$.
It is proved in \cite{LM13} (Lemma 15) that $S_t^0$ is an interval for the Bruhat order of $S_t$,
{\sl i.e.} there exists
an element  $w'_0\in S_t$ such that $S_t^0=\{w \in S_t \, \vert\, w\leq w'_0\}$.
Furthermore, $w'_0$ avoids $312$ pattern, {\sl i.e.} there does not exist a triple $i<j<k$ such that
$w'_0(k)<w'_0(i)<w'_0(j)$ and  $\caI_{w'_0}$ is irreducible.

We denote by $\caM$ the multisegment  $([a_1; b_1], \ldots , [a_t; b_t ])$ and for any
$w \in S_t^0$, we set $\caM_w=  ([a_{w(1)}; b_1], \ldots  [a_{w(t)}; b_t ])$.

\begin{prop}\label{prop:zelorder} The map $w\mapsto\caM_w$ is a bijection from $S_t^0$ onto the set
of multisegments $\frn$ with $\frn\leq\caM$. If $w_1\leq w_2$ in $S_t^0$, then $\caM_{w_2}\leq\caM_{w_1}$.
\end{prop}

\begin{proof}
Here the empty segments of $\caM_w$ are omitted. Since the $a_j$ are pairwise distinct, the pairing
of the ends $b_i$ with the beginnings $a_{w(i)}$ can be recovered from $\caM_w$: an end $b_i$ which
is not the end of a segment of $\caM_w$ is paired with the beginning $a_j=b_i+1$. Hence the map is injective.

Let $w_1<w_2$ in $S_t^0$. We want to show that $\caM_{w_2}\leq\caM_{w_1}$. Since the order on
multisegments is transitive, and since, $S_t^0$ being a lower interval, $w_1$ and $w_2$ are joined by
a chain of elements of $S_t^0$ in which each element covers the previous one (\cite{BB05}, Theorem
2.2.6), we may assume that $w_2$ covers $w_1$. Then $w_2=w_1\cdot(i,j)$ with $i<j$ and
$w_1(i)<w_1(j)$ (\cite{BB05}, Lemma 2.1.4). Then $a_{w_1(i)}>a_{w_1(j)}$, $b_i>b_j$, and $a_{w_1(i)}\leq b_j+1$ since $w_2\in
S_t^0$. Hence the segments $[a_{w_1(i)};b_i]$ and $[a_{w_1(j)};b_j]$ are linked, and $\caM_{w_2}$
is obtained from $\caM_{w_1}$ by replacing them by their union $[a_{w_1(j)};b_i]$ and their
intersection $[a_{w_1(i)};b_j]$ (empty if $a_{w_1(i)}=b_j+1$). Thus $\caM_{w_2}\leq\caM_{w_1}$.
Taking $w_1=\Id$, we get $\caM_w\leq\caM$ for every $w\in S_t^0$.

Conversely, it suffices to show that a multisegment $\frn$ obtained from some $\caM_w$, $w\in
S_t^0$, by one elementary operation is of the form $\caM_x$, $x\in S_t^0$. Let $[a;b]$, $[a';b']$
be the two linked segments, with $a<a'\leq b+1\leq b'$. They are replaced by $[a;b']$ and $[a';b]$,
the latter being empty if $a'=b+1$, in which case we keep it formally as the empty segment
$[b+1;b]$. The multisets of beginnings and of ends are preserved. Hence $\frn$ pairs bijectively
the ends $b_i$ with the beginnings $a_j$, with $a_j\leq b_i+1$. That is, $\frn=\caM_x$ with $x\in S_t^0$.
\end{proof}

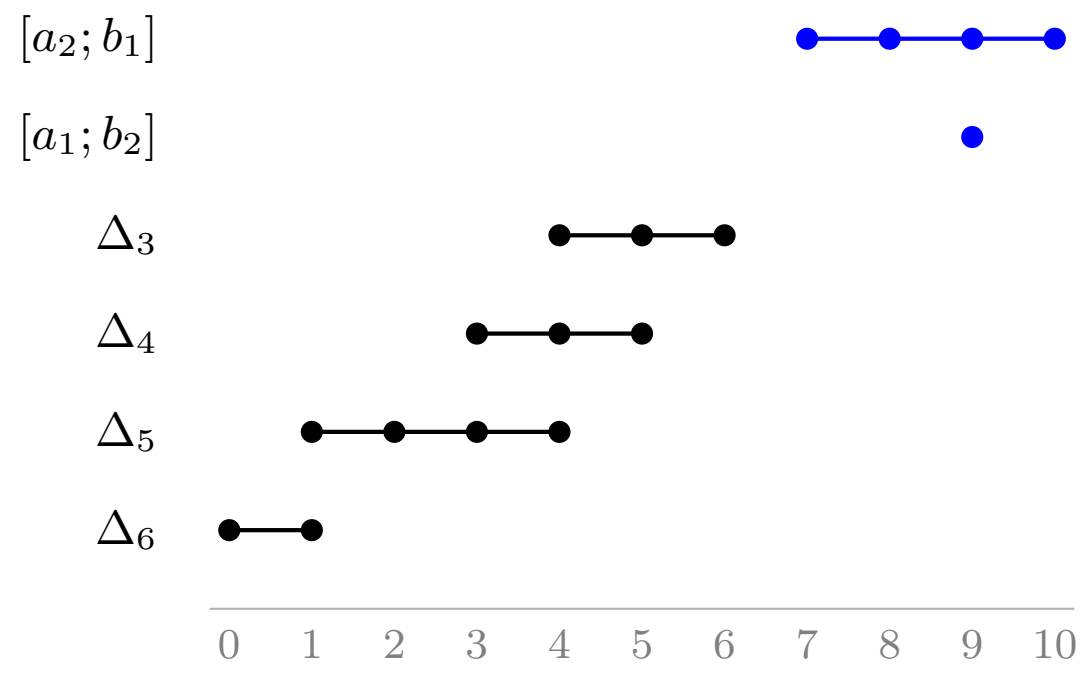
\begin{figure}[h]
\centering
\begin{tikzpicture}[scale=1.3]
\draw[blue,thick] (2.94,-0.50)--(4.20,-0.50);
\fill[blue] (2.94,-0.50) circle (1.6pt);
\fill[blue] (3.36,-0.50) circle (1.6pt);
\fill[blue] (3.78,-0.50) circle (1.6pt);
\fill[blue] (4.20,-0.50) circle (1.6pt);
\node[left,font=\scriptsize] at (-0.25,-0.50) {$[a_2;b_1]$};
\draw[blue,thick] (3.78,-1.00)--(3.78,-1.00);
\fill[blue] (3.78,-1.00) circle (1.6pt);
\node[left,font=\scriptsize] at (-0.25,-1.00) {$[a_1;b_2]$};
\draw[black,thick] (1.68,-1.50)--(2.52,-1.50);
\fill[black] (1.68,-1.50) circle (1.6pt);
\fill[black] (2.10,-1.50) circle (1.6pt);
\fill[black] (2.52,-1.50) circle (1.6pt);
\node[left,font=\scriptsize] at (-0.25,-1.50) {$\Delta_{3}$};
\draw[black,thick] (1.26,-2.00)--(2.10,-2.00);
\fill[black] (1.26,-2.00) circle (1.6pt);
\fill[black] (1.68,-2.00) circle (1.6pt);
\fill[black] (2.10,-2.00) circle (1.6pt);
\node[left,font=\scriptsize] at (-0.25,-2.00) {$\Delta_{4}$};
\draw[black,thick] (0.42,-2.50)--(1.68,-2.50);
\fill[black] (0.42,-2.50) circle (1.6pt);
\fill[black] (0.84,-2.50) circle (1.6pt);
\fill[black] (1.26,-2.50) circle (1.6pt);
\fill[black] (1.68,-2.50) circle (1.6pt);
\node[left,font=\scriptsize] at (-0.25,-2.50) {$\Delta_{5}$};
\draw[black,thick] (0.00,-3.00)--(0.42,-3.00);
\fill[black] (0.00,-3.00) circle (1.6pt);
\fill[black] (0.42,-3.00) circle (1.6pt);
\node[left,font=\scriptsize] at (-0.25,-3.00) {$\Delta_{6}$};
\draw[gray!60] (-0.10,-3.40)--(4.30,-3.40);
\node[font=\tiny,gray] at (0.00,-3.58) {0};
\node[font=\tiny,gray] at (0.42,-3.58) {1};
\node[font=\tiny,gray] at (0.84,-3.58) {2};
\node[font=\tiny,gray] at (1.26,-3.58) {3};
\node[font=\tiny,gray] at (1.68,-3.58) {4};
\node[font=\tiny,gray] at (2.10,-3.58) {5};
\node[font=\tiny,gray] at (2.52,-3.58) {6};
\node[font=\tiny,gray] at (2.94,-3.58) {7};
\node[font=\tiny,gray] at (3.36,-3.58) {8};
\node[font=\tiny,gray] at (3.78,-3.58) {9};
\node[font=\tiny,gray] at (4.20,-3.58) {10};
\end{tikzpicture}
\caption{The multisegment $\caM_{s_1}$ of Example \ref{ex:run2}; the two modified segments are in
blue.}\label{fig:cover}
\end{figure}

\begin{example}\label{ex:run2} For the ladder of Example \ref{ex:run}, $(b_i+1)_i=(11,10,7,6,5,2)$,
and the conditions $a_{w(i)}\leq b_i+1$ defining $S_6^0$ read
\[w(3)\geq2,\qquad w(4)\geq3,\qquad w(5)\geq3,\qquad w(6)\geq5.\]
The set $S_6^0$ has $48$ elements. It is the interval $[\Id,w'_0]$ with $w'_0=612435$ (one-line
notation), of length $6$, and its elements of length $0,1,\ldots,6$ number $1,5,11,14,11,5,1$.
Theorem \ref{BGG} thus predicts an exact sequence
\[0\rightarrow\pi\rightarrow\caJ_{\Id}\rightarrow\caJ^{\oplus}_{(5)}\rightarrow\caJ^{\oplus}_{(11)}
\rightarrow\caJ^{\oplus}_{(14)}\rightarrow\caJ^{\oplus}_{(11)}\rightarrow\caJ^{\oplus}_{(5)}\rightarrow
\caJ_{w'_0}\rightarrow0,\]
where $\caJ^\oplus_{(n)}$ denotes a sum of $n$ costandard modules. The five elements of length $1$
are $s_1,\ldots,s_5$. For instance, $\caM_{s_1}$ is obtained from $\caM=\caM_{\Id}$ by the
elementary operation on the linked segments $[9;10]$ and $[7;9]$, which are replaced by their union
$[7;10]$ and their intersection $[9;9]$ (Figure \ref{fig:cover}).
\end{example}

\begin{lemma}\label{lem:twoseg} Let $\beta<\alpha\leq b'+1\leq b$, and put
$\Delta=\Delta([\alpha;b])$, $\Delta'=\Delta([\beta;b'])$, $\Delta_\cup=\Delta([\beta;b])$ and
$\Delta_\cap=\Delta([\alpha;b'])$, with $\Delta_\cap=1$ if $\alpha=b'+1$. Then
$\Delta_\cup\times\Delta_\cap\simeq\Delta_\cap\times\Delta_\cup$ is irreducible, and there are exact sequences
\[0\rightarrow\Delta_\cup\times\Delta_\cap\rightarrow\Delta\times\Delta'\rightarrow
L(\Delta,\Delta')\rightarrow0,\qquad 0\rightarrow
L(\Delta,\Delta')\rightarrow\Delta'\times\Delta\rightarrow\Delta_\cup\times\Delta_\cap\rightarrow0.\]
\end{lemma}

\begin{proof}
This is well known: it is the analogue, for essentially square-integrable representations, of
\cite{Zel80}, Proposition 4.6, from which it follows by the Zelevinsky involution (\cite{Zel80},
9.12--9.15, which exchanges the representations $\langle\Delta\rangle$ of \emph{loc. cit.} and the
$\Delta(\Delta)$, together with \cite{Aub95}, \cite{MW86}, by which it maps irreducible representations
to irreducible representations): $\Delta\times\Delta'$ and $\Delta'\times\Delta$ have length two, with
composition factors $L(\Delta,\Delta')$ and $\Delta_\cup\times\Delta_\cap$, the latter being irreducible
by (Z). The exact sequences follow, since $L(\Delta,\Delta')$ is the unique irreducible quotient of
$\Delta\times\Delta'$ \cite{Rod82} and the unique irreducible subrepresentation of $\Delta'\times\Delta$
by (L).
\end{proof}

\begin{prop}\label{prop:inclusion}
 If $w_1\leq w_2$ are elements in $S_t^0$, then $\caI_{w_2}$ is a subrepresentation of $\caI_{w_1}$
 \end{prop}

\proof By the chain property of the Bruhat order (\cite{BB05}, Theorem 2.2.6), applied in the
interval $S_t^0$,
 it suffices to treat the case of a covering $w_1\lessdot w_2$, i.e.\ $w_1<w_2$ and there is no
 $w$ with $w_1<w<w_2$; equivalently $w_1<w_2$ and $l(w_2)=l(w_1)+1$, the Bruhat order being graded
 by the length (\cite{BB05}, \S2.2).
  By \cite{BB05}, Lemma 2.1.4, we then have
  $w_2=w_1\cdot(i,j)$ with $i<j$, $w_1(i)<w_1(j)$, and for every $k$ with $i<k<j$,
   either $w_1(k)<w_1(i)$ or $w_1(k)>w_1(j)$. Write $\Delta^{w_1}_k=\Delta([a_{w_1(k)};b_k])$, and let
\[\Delta_\cup=\Delta([a_{w_1(j)};b_i]),\qquad \Delta_\cap=\Delta([a_{w_1(i)};b_j])\]
be the segments of $\caM_{w_2}$ at the positions $i$ and $j$. Let $K_1$ (resp. $K_2$) be the set of
$k\in\{i+1,\ldots,j-1\}$ with $w_1(k)<w_1(i)$ (resp. $w_1(k)>w_1(j)$). Recall that $a_1>\cdots>a_t$
and $b_1>\cdots>b_t$.
\begin{itemize}
\item If $k\in K_1$, then $a_{w_1(k)}>a_{w_1(i)}>a_{w_1(j)}$ and $b_k<b_i$. Hence the segment of
$\Delta^{w_1}_k$ is contained in those of $\Delta^{w_1}_i$ and $\Delta_\cup$.
\item If $k\in K_2$, then $a_{w_1(k)}<a_{w_1(j)}<a_{w_1(i)}$ and $b_k>b_j$. Hence the segment of
$\Delta^{w_1}_k$ contains those of $\Delta^{w_1}_j$ and $\Delta_\cap$.
\item If $k_2<k_1$ with $k_2\in K_2$ and $k_1\in K_1$, then the segment of $\Delta^{w_1}_{k_1}$
is contained in that of $\Delta^{w_1}_{k_2}$.
\end{itemize}
A segment contained in another one is not linked with it, so the corresponding product is irreducible
and the two factors commute by (Z). Using only the commutations listed above, we obtain
\[\caI_{w_1}\simeq X\times\Big(\prod_{k\in
K_1}\Delta^{w_1}_k\Big)\times\Delta^{w_1}_i\times\Delta^{w_1}_j\times
\Big(\prod_{k\in K_2}\Delta^{w_1}_k\Big)\times Y,\]
where $X$ (resp. $Y$) is the product of the factors of $\caI_{w_1}$ of index $<i$ (resp. $>j$),
and the products over $K_1$ and $K_2$ are taken in increasing order of $k$. Since $w_2\in S_t^0$,
 we have $a_{w_1(i)}\leq b_j+1$. Hence the segments of $\Delta^{w_1}_i$ and $\Delta^{w_1}_j$ are
  linked or juxtaposed, $\Delta_\cup$ and $\Delta_\cap$ are their union and intersection, and
  $\Delta_\cup\times\Delta_\cap\hookrightarrow\Delta^{w_1}_i\times\Delta^{w_1}_j$ by Lemma
  \ref{lem:twoseg}, with
  $\Delta_\cap=1$ in the juxtaposed case. Parabolic induction is exact, so we get an embedding
\[X\times\Big(\prod_{K_1}\Delta^{w_1}_k\Big)\times\Delta_\cup\times\Delta_\cap\times
\Big(\prod_{K_2}\Delta^{w_1}_k\Big)\times Y\hookrightarrow\caI_{w_1}.\]
By the same inclusions of segments, the left-hand side is isomorphic to $\caI_{w_2}$:
the factors indexed by $K_1$ move to the right of $\Delta_\cup$, those indexed by $K_2$ move
to the left of $\Delta_\cap$, and the original relative order is then restored. This proves the proposition.
 Reversing the order of all the factors, the same argument gives a surjection
 $\caJ_{w_1}\twoheadrightarrow\caJ_{w_2}$, using the surjection
 $\Delta^{w_1}_j\times\Delta^{w_1}_i\twoheadrightarrow\Delta_\cup\times\Delta_\cap\simeq\Delta_\cap\times\Delta_\cup$
 of Lemma \ref{lem:twoseg}.
\qed

\section{Minors, words and multiplicities}\label{sec:reduction}

Throughout this section and the next one, $\pi=L(\Delta_1,\ldots,\Delta_t)$,
$\Delta_i=\Delta([a_i;b_i])$, is a ladder representation, with $a_1>\cdots>a_t$ and $b_1>\cdots>b_t$.
This section contains the combinatorial and representation-theoretic results on which the proof of
Theorem \ref{BGG} rests; the proof itself, by homological algebra, is in \S\ref{sec:complex}. We prove
Theorem \ref{BGG} by induction on $t$, the case
 $t=1$ being trivial ($\pi=\Delta_1=\caJ_{\Id}$). The strategy mimics the Cousin resolution: the
 interval $S_t^0$
 is stratified according to the value $w(1)$, the exactness of parabolic induction plays the role
 of the exactness
  of direct images along affine embeddings, and the key lemma (Proposition \ref{conj:key}) plays the role
  of the purity of local cohomology. The two representation-theoretic inputs, the key lemma and the
  multiplicity one statement of Proposition \ref{prop:multone}, are proved by counting words in
  minimal Jacquet modules (\S\S\ref{ssec:words}-\ref{ssec:multone}).

\subsection{Costandard modules}
For $w\in S_t^0$, recall that
\[\caJ_w=\Delta([a_{w(t)};b_t])\times\cdots\times\Delta([a_{w(1)};b_1]),\]
where empty segments are omitted. The segments of $\caJ_w$ are written in increasing order
of their ends; hence $\caJ_w^\vee$ is a standard module and $\caJ_w$ has a unique irreducible
subrepresentation,
namely $L(\caM_w)$. More generally, if $\sigma=L(\Delta''_1,\ldots,\Delta''_r)$ is a ladder
representation, $\sigma$
is the unique irreducible subrepresentation of $\Delta''_r\times\cdots\times\Delta''_1$.

\subsection{Generalized ladders}\label{ssec:generalized}
We say that $(\Delta_1,\ldots,\Delta_t)$, with $\Delta_i=\Delta([a_i;b_i])$, is a generalized ladder
if $a_1>\cdots>a_t$, $b_1>\cdots>b_t$ and $a_i\leq b_i+1$ for all $i$. Empty segments $[a_i;a_i-1]$
 are thus allowed. We define $S_t^0$, $\caM_w$, $\caI_w$, $\caJ_w$ by the same formulas as before,
 empty segments being omitted in the products. We also set $L(\Delta_1,\ldots,\Delta_t)=L(\caM)$,
 empty segments being omitted.

\begin{lemma}\label{lem:empty} Let $(\Delta_1,\ldots,\Delta_t)$ be a generalized ladder with $a_i=b_i+1$
 for some $i$. Then every $w\in S_t^0$ satisfies $w(i)=i$ and $w(\{1,\ldots,i-1\})=\{1,\ldots,i-1\}$.
 Removing the index $i$ gives a bijection $w\mapsto w^\flat$ from $S_t^0$ onto the set $S^0_{t-1}$
 attached to $(\Delta_1,\ldots,\widehat{\Delta_i},\ldots,\Delta_t)$. This bijection preserves
 length and Bruhat order,
 and $\caI_w=\caI_{w^\flat}$, $\caJ_w=\caJ_{w^\flat}$. Consequently, Theorem \ref{LM} (ii) and
 Theorem \ref{BGG} hold for generalized ladders as soon as they hold for ladders with fewer segments.
\end{lemma}

\begin{proof}
Let $w\in S_t^0$. For $i'\geq i$ we have $a_{w(i')}\leq b_{i'}+1\leq b_i+1=a_i$, hence $w(i')\geq i$.
Thus $w$ preserves $\{i,\ldots,t\}$, hence also $\{1,\ldots,i-1\}$. If $w(i)\neq i$, then $i=w(i')$
for some $i'>i$,
and $a_i\leq b_{i'}+1<b_i+1=a_i$, a contradiction. Therefore $w$ lies in the standard parabolic subgroup
$S_{i-1}\times S_1\times S_{t-i}$, on which length and Bruhat order are those of the product. Conversely,
inserting the fixed point $i$ into an element of $S^0_{t-1}$ gives an element of $S^0_t$, since the condition
 at $i$ is $a_i\leq b_i+1$. The segment indexed by $i$ in $\caM_w$ is the empty segment $[a_i;b_i]$,
 which gives the last assertions.
\end{proof}

By Lemma \ref{lem:empty}, the results of this section, stated for ladders, hold for generalized
ladders as well. We use this for the minors defined below.

\subsection{Minors}\label{ssec:minors}
For $k=1,\ldots,t$, set $\Delta'_k=\Delta([a_k;b_1])$ and, for $i=2,\ldots,t$,
\[\Delta^{(k)}_i=\begin{cases}\Delta([a_{i-1};b_i]) & \text{if } 2\leq i\leq k,\\ \Delta([a_i;b_i])
& \text{if } k<i\leq t.\end{cases}\]
The terminology comes from Tadi\'c's formula. In the commutative ring $\caR$, Theorem \ref{LM} (ii)
says that $[\pi]$ is the determinant of the matrix $(\Delta([a_j;b_i]))_{1\leq i,j\leq t}$, with the
conventions of the introduction for $b_i<a_j$. Its expansion along the first row is
$\sum_k(-1)^{k-1}[\Delta'_k]\,D_k$, where $D_k$ is the minor obtained by deleting the first row and the
column $k$; this minor is the determinant attached in the same way to the segments
$\Delta^{(k)}_2,\ldots,\Delta^{(k)}_t$ (beginnings $a_1,\ldots,\widehat{a_k},\ldots,a_t$, ends
$b_2,\ldots,b_t$). We say that the minor of index $k$ vanishes if $a_{i-1}>b_i+1$ for some $2\leq i\leq k$, and
we set $\pi_k=0$ in this case. Otherwise, $(\Delta^{(k)}_2,\ldots,\Delta^{(k)}_t)$ is a generalized
 ladder (\S\ref{ssec:generalized}) and we set $\pi_k=L(\Delta^{(k)}_2,\ldots,\Delta^{(k)}_t)$. Let
 $c_k\in S_t$ be defined by
\[c_k(1)=k,\qquad c_k(i)=i-1\ (2\leq i\leq k),\qquad c_k(i)=i\ (k<i\leq t),\]
so that $c_1=\Id$, $l(c_k)=k-1$, and $\caM_{c_k}=(\Delta'_k,\Delta^{(k)}_2,\ldots,\Delta^{(k)}_t)$.

\begin{lemma}\label{lem:minors} Let $1\leq k\leq t$.
\begin{enumerate}
\item[(i)] If $\pi_k\neq 0$, then $(\Delta^{(k)}_2,\ldots,\Delta^{(k)}_t)$ is a generalized ladder
with $t-1$ segments, and $a_k\leq b_1$, so that $\Delta'_k\neq 0$.
\item[(ii)] If the minor of index $k$ vanishes, so does the minor of index $k'$ for every $k'\geq k$.
\item[(iii)] Let $S_{t,k}^0=\{w\in S_t^0 : w(1)=k\}$. Then $S_{t,k}^0\neq\emptyset$ if and only if
$\pi_k\neq 0$, if and only if $c_k\in S_t^0$. In that case, $c_k$ is the unique element of minimal
length of $S_{t,k}^0$.
\item[(iv)] Assume $\pi_k\neq 0$. For $w\in S_{t,k}^0$, let $\bar w\in S_{t-1}$ be the
standardization of the sequence $(w(2),\allowbreak\ldots,\allowbreak w(t))$, i.e.\ the unique permutation $\bar w$ of
$\{1,\ldots,t-1\}$ such that $\bar w(i)<\bar w(j)$ if and only if $w(i+1)<w(j+1)$.
Then $w\mapsto\bar w$ is a bijection from $S_{t,k}^0$ onto the set $S^0_{t-1}$ attached to the
generalized ladder
$(\Delta^{(k)}_2,\ldots,\Delta^{(k)}_t)$, it is compatible with the Bruhat orders, $l(w)=(k-1)+l(\bar w)$, and
\[\caJ_w\simeq \caJ^{(k)}_{\bar w}\times\Delta'_k,\]
where $\caJ^{(k)}$ denotes the costandard modules attached to $\pi_k$.
\item[(v)] If $w\leq w'$ in $S_t$, then $w(1)\leq w'(1)$.
\item[(vi)] If $c_{k+1}\in S_t^0$, then $c_k\lessdot c_{k+1}$. If $w\in S_t^0$ satisfies $w(1)\geq
k$, then $l(w)\geq k-1$, with equality if and only if $w=c_k$.
\end{enumerate}
\end{lemma}

\begin{proof}
(i) The beginnings $a_1>\cdots>a_{k-1}>a_{k+1}>\cdots>a_t$ and the ends $b_2>\cdots>b_t$ of the segments
$\Delta^{(k)}_i$ are strictly decreasing, and $a_k\leq a_1\leq b_1$. (ii) is clear from the
definition. (iii) If $a_{i-1}>b_i+1$
for some $2\leq i\leq k$, then $a_j>b_{i'}+1$ for all $j\leq i-1$ and $i'\geq i$. An element $w\in
S^0_{t,k}$ would then
send the $t-i+1$ indices $i,\ldots,t$ into $\{i,\ldots,t\}\setminus\{k\}$, which has $t-i$
elements: this is impossible.
Conversely, if the minor does not vanish, then $c_k\in S_t^0$ by definition. The other assertions
of (iii) and (iv)
follow from the identity $l(w)=\#\{j>1 : w(j)<w(1)\}+l(\bar w)$ and from the fact that the
condition $a_{w(i)}\leq b_i+1$
 for $i\geq 2$ is exactly the condition defining the interval attached to $\pi_k$. The isomorphism
 in (iv) holds because
 $\Delta([a_{w(1)};b_1])=\Delta'_k$ is the rightmost factor of $\caJ_w$. The compatibility with the
 Bruhat orders in (iv) also
 follows from the tableau criterion (\cite{BB05}, Theorem 2.6.3 and Proposition 2.4.8), in the
 form: $x\leq y$ if and only if,
 for every $j$, the set $x(\{1,\ldots,j\})$ is smaller than $y(\{1,\ldots,j\})$ in the Gale order.
 Recall that $A\leq B$ in the
 Gale order if and only if $\#\{a\in A : a\geq m\}\leq\#\{b\in B : b\geq m\}$ for every $m$. If
 $x(1)=y(1)=k$, the initial
 sets of $x$ and $y$ are obtained from those of $\bar x$ and $\bar y$ by adding the common element $k$,
 after an increasing relabelling of the values. Both operations preserve and reflect the Gale
 order. (v) is the case $j=1$ of the same criterion. (vi) We have $c_{k+1}=c_k\cdot(1,k+1)$ and
 $l(c_{k+1})=l(c_k)+1$. If $w(1)=k'\geq k$, then $l(w)=(k'-1)+l(\bar w)\geq k-1$, with equality if
 and only if $k'=k$ and $\bar w=\Id$, i.e. $w=c_k$.
\end{proof}

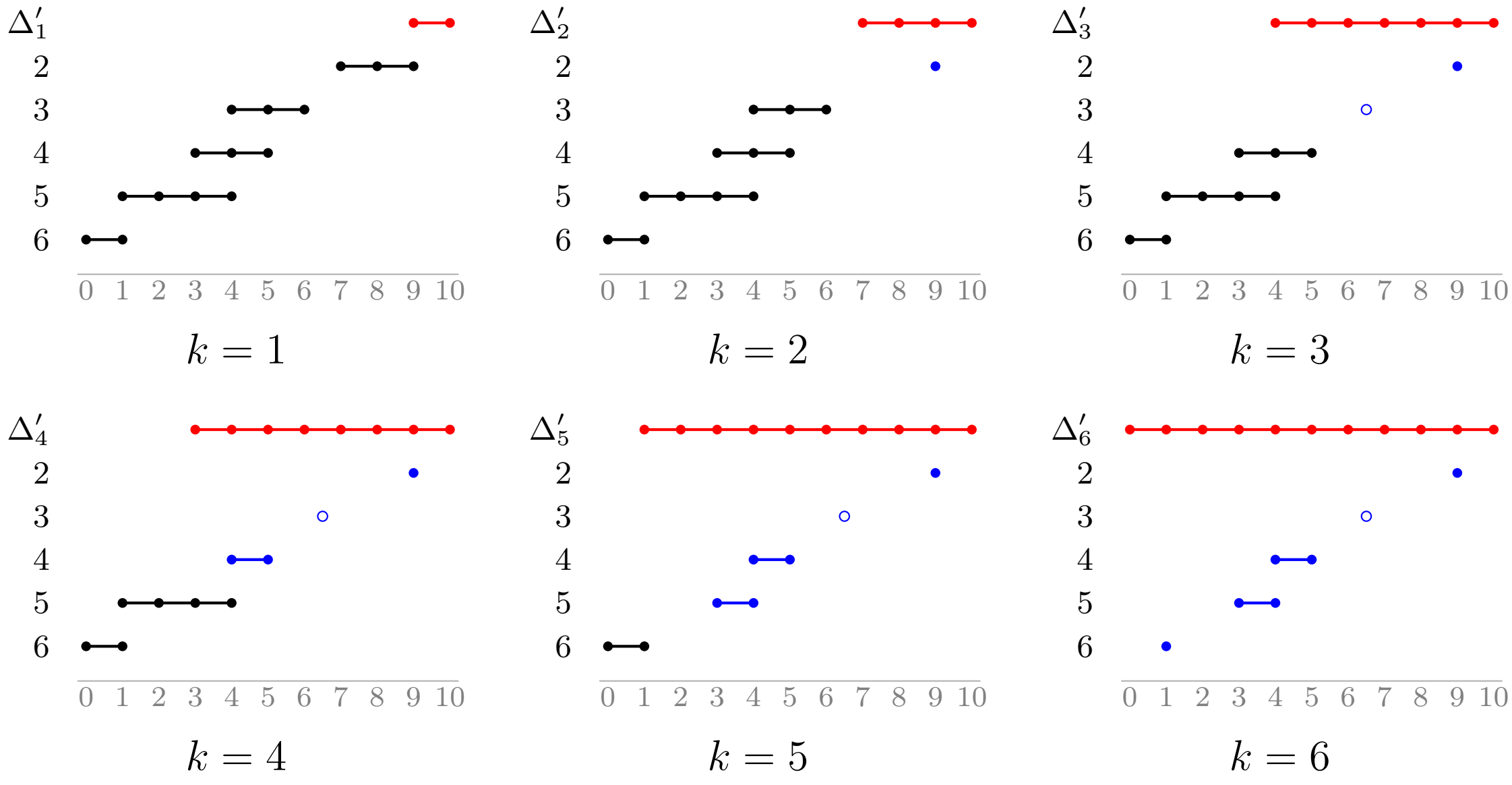
\begin{figure}[h]
\centering
\begin{tabular}{ccc}
\begin{tikzpicture}[scale=0.85]
\draw[red,thick] (3.78,-0.50)--(4.20,-0.50);
\fill[red] (3.78,-0.50) circle (1.6pt);
\fill[red] (4.20,-0.50) circle (1.6pt);
\node[left,font=\scriptsize] at (-0.25,-0.50) {$\Delta'_{1}$};
\draw[black,thick] (2.94,-1.00)--(3.78,-1.00);
\fill[black] (2.94,-1.00) circle (1.6pt);
\fill[black] (3.36,-1.00) circle (1.6pt);
\fill[black] (3.78,-1.00) circle (1.6pt);
\node[left,font=\scriptsize] at (-0.25,-1.00) {$2$};
\draw[black,thick] (1.68,-1.50)--(2.52,-1.50);
\fill[black] (1.68,-1.50) circle (1.6pt);
\fill[black] (2.10,-1.50) circle (1.6pt);
\fill[black] (2.52,-1.50) circle (1.6pt);
\node[left,font=\scriptsize] at (-0.25,-1.50) {$3$};
\draw[black,thick] (1.26,-2.00)--(2.10,-2.00);
\fill[black] (1.26,-2.00) circle (1.6pt);
\fill[black] (1.68,-2.00) circle (1.6pt);
\fill[black] (2.10,-2.00) circle (1.6pt);
\node[left,font=\scriptsize] at (-0.25,-2.00) {$4$};
\draw[black,thick] (0.42,-2.50)--(1.68,-2.50);
\fill[black] (0.42,-2.50) circle (1.6pt);
\fill[black] (0.84,-2.50) circle (1.6pt);
\fill[black] (1.26,-2.50) circle (1.6pt);
\fill[black] (1.68,-2.50) circle (1.6pt);
\node[left,font=\scriptsize] at (-0.25,-2.50) {$5$};
\draw[black,thick] (0.00,-3.00)--(0.42,-3.00);
\fill[black] (0.00,-3.00) circle (1.6pt);
\fill[black] (0.42,-3.00) circle (1.6pt);
\node[left,font=\scriptsize] at (-0.25,-3.00) {$6$};
\draw[gray!60] (-0.10,-3.40)--(4.30,-3.40);
\node[font=\tiny,gray] at (0.00,-3.58) {0};
\node[font=\tiny,gray] at (0.42,-3.58) {1};
\node[font=\tiny,gray] at (0.84,-3.58) {2};
\node[font=\tiny,gray] at (1.26,-3.58) {3};
\node[font=\tiny,gray] at (1.68,-3.58) {4};
\node[font=\tiny,gray] at (2.10,-3.58) {5};
\node[font=\tiny,gray] at (2.52,-3.58) {6};
\node[font=\tiny,gray] at (2.94,-3.58) {7};
\node[font=\tiny,gray] at (3.36,-3.58) {8};
\node[font=\tiny,gray] at (3.78,-3.58) {9};
\node[font=\tiny,gray] at (4.20,-3.58) {10};
\end{tikzpicture}&
\begin{tikzpicture}[scale=0.85]
\draw[red,thick] (2.94,-0.50)--(4.20,-0.50);
\fill[red] (2.94,-0.50) circle (1.6pt);
\fill[red] (3.36,-0.50) circle (1.6pt);
\fill[red] (3.78,-0.50) circle (1.6pt);
\fill[red] (4.20,-0.50) circle (1.6pt);
\node[left,font=\scriptsize] at (-0.25,-0.50) {$\Delta'_{2}$};
\draw[blue,thick] (3.78,-1.00)--(3.78,-1.00);
\fill[blue] (3.78,-1.00) circle (1.6pt);
\node[left,font=\scriptsize] at (-0.25,-1.00) {$2$};
\draw[black,thick] (1.68,-1.50)--(2.52,-1.50);
\fill[black] (1.68,-1.50) circle (1.6pt);
\fill[black] (2.10,-1.50) circle (1.6pt);
\fill[black] (2.52,-1.50) circle (1.6pt);
\node[left,font=\scriptsize] at (-0.25,-1.50) {$3$};
\draw[black,thick] (1.26,-2.00)--(2.10,-2.00);
\fill[black] (1.26,-2.00) circle (1.6pt);
\fill[black] (1.68,-2.00) circle (1.6pt);
\fill[black] (2.10,-2.00) circle (1.6pt);
\node[left,font=\scriptsize] at (-0.25,-2.00) {$4$};
\draw[black,thick] (0.42,-2.50)--(1.68,-2.50);
\fill[black] (0.42,-2.50) circle (1.6pt);
\fill[black] (0.84,-2.50) circle (1.6pt);
\fill[black] (1.26,-2.50) circle (1.6pt);
\fill[black] (1.68,-2.50) circle (1.6pt);
\node[left,font=\scriptsize] at (-0.25,-2.50) {$5$};
\draw[black,thick] (0.00,-3.00)--(0.42,-3.00);
\fill[black] (0.00,-3.00) circle (1.6pt);
\fill[black] (0.42,-3.00) circle (1.6pt);
\node[left,font=\scriptsize] at (-0.25,-3.00) {$6$};
\draw[gray!60] (-0.10,-3.40)--(4.30,-3.40);
\node[font=\tiny,gray] at (0.00,-3.58) {0};
\node[font=\tiny,gray] at (0.42,-3.58) {1};
\node[font=\tiny,gray] at (0.84,-3.58) {2};
\node[font=\tiny,gray] at (1.26,-3.58) {3};
\node[font=\tiny,gray] at (1.68,-3.58) {4};
\node[font=\tiny,gray] at (2.10,-3.58) {5};
\node[font=\tiny,gray] at (2.52,-3.58) {6};
\node[font=\tiny,gray] at (2.94,-3.58) {7};
\node[font=\tiny,gray] at (3.36,-3.58) {8};
\node[font=\tiny,gray] at (3.78,-3.58) {9};
\node[font=\tiny,gray] at (4.20,-3.58) {10};
\end{tikzpicture}&
\begin{tikzpicture}[scale=0.85]
\draw[red,thick] (1.68,-0.50)--(4.20,-0.50);
\fill[red] (1.68,-0.50) circle (1.6pt);
\fill[red] (2.10,-0.50) circle (1.6pt);
\fill[red] (2.52,-0.50) circle (1.6pt);
\fill[red] (2.94,-0.50) circle (1.6pt);
\fill[red] (3.36,-0.50) circle (1.6pt);
\fill[red] (3.78,-0.50) circle (1.6pt);
\fill[red] (4.20,-0.50) circle (1.6pt);
\node[left,font=\scriptsize] at (-0.25,-0.50) {$\Delta'_{3}$};
\draw[blue,thick] (3.78,-1.00)--(3.78,-1.00);
\fill[blue] (3.78,-1.00) circle (1.6pt);
\node[left,font=\scriptsize] at (-0.25,-1.00) {$2$};
\draw[blue] (2.73,-1.50) circle (1.6pt);
\node[left,font=\scriptsize] at (-0.25,-1.50) {$3$};
\draw[black,thick] (1.26,-2.00)--(2.10,-2.00);
\fill[black] (1.26,-2.00) circle (1.6pt);
\fill[black] (1.68,-2.00) circle (1.6pt);
\fill[black] (2.10,-2.00) circle (1.6pt);
\node[left,font=\scriptsize] at (-0.25,-2.00) {$4$};
\draw[black,thick] (0.42,-2.50)--(1.68,-2.50);
\fill[black] (0.42,-2.50) circle (1.6pt);
\fill[black] (0.84,-2.50) circle (1.6pt);
\fill[black] (1.26,-2.50) circle (1.6pt);
\fill[black] (1.68,-2.50) circle (1.6pt);
\node[left,font=\scriptsize] at (-0.25,-2.50) {$5$};
\draw[black,thick] (0.00,-3.00)--(0.42,-3.00);
\fill[black] (0.00,-3.00) circle (1.6pt);
\fill[black] (0.42,-3.00) circle (1.6pt);
\node[left,font=\scriptsize] at (-0.25,-3.00) {$6$};
\draw[gray!60] (-0.10,-3.40)--(4.30,-3.40);
\node[font=\tiny,gray] at (0.00,-3.58) {0};
\node[font=\tiny,gray] at (0.42,-3.58) {1};
\node[font=\tiny,gray] at (0.84,-3.58) {2};
\node[font=\tiny,gray] at (1.26,-3.58) {3};
\node[font=\tiny,gray] at (1.68,-3.58) {4};
\node[font=\tiny,gray] at (2.10,-3.58) {5};
\node[font=\tiny,gray] at (2.52,-3.58) {6};
\node[font=\tiny,gray] at (2.94,-3.58) {7};
\node[font=\tiny,gray] at (3.36,-3.58) {8};
\node[font=\tiny,gray] at (3.78,-3.58) {9};
\node[font=\tiny,gray] at (4.20,-3.58) {10};
\end{tikzpicture}\\
$k=1$ & $k=2$ & $k=3$\\[6pt]
\begin{tikzpicture}[scale=0.85]
\draw[red,thick] (1.26,-0.50)--(4.20,-0.50);
\fill[red] (1.26,-0.50) circle (1.6pt);
\fill[red] (1.68,-0.50) circle (1.6pt);
\fill[red] (2.10,-0.50) circle (1.6pt);
\fill[red] (2.52,-0.50) circle (1.6pt);
\fill[red] (2.94,-0.50) circle (1.6pt);
\fill[red] (3.36,-0.50) circle (1.6pt);
\fill[red] (3.78,-0.50) circle (1.6pt);
\fill[red] (4.20,-0.50) circle (1.6pt);
\node[left,font=\scriptsize] at (-0.25,-0.50) {$\Delta'_{4}$};
\draw[blue,thick] (3.78,-1.00)--(3.78,-1.00);
\fill[blue] (3.78,-1.00) circle (1.6pt);
\node[left,font=\scriptsize] at (-0.25,-1.00) {$2$};
\draw[blue] (2.73,-1.50) circle (1.6pt);
\node[left,font=\scriptsize] at (-0.25,-1.50) {$3$};
\draw[blue,thick] (1.68,-2.00)--(2.10,-2.00);
\fill[blue] (1.68,-2.00) circle (1.6pt);
\fill[blue] (2.10,-2.00) circle (1.6pt);
\node[left,font=\scriptsize] at (-0.25,-2.00) {$4$};
\draw[black,thick] (0.42,-2.50)--(1.68,-2.50);
\fill[black] (0.42,-2.50) circle (1.6pt);
\fill[black] (0.84,-2.50) circle (1.6pt);
\fill[black] (1.26,-2.50) circle (1.6pt);
\fill[black] (1.68,-2.50) circle (1.6pt);
\node[left,font=\scriptsize] at (-0.25,-2.50) {$5$};
\draw[black,thick] (0.00,-3.00)--(0.42,-3.00);
\fill[black] (0.00,-3.00) circle (1.6pt);
\fill[black] (0.42,-3.00) circle (1.6pt);
\node[left,font=\scriptsize] at (-0.25,-3.00) {$6$};
\draw[gray!60] (-0.10,-3.40)--(4.30,-3.40);
\node[font=\tiny,gray] at (0.00,-3.58) {0};
\node[font=\tiny,gray] at (0.42,-3.58) {1};
\node[font=\tiny,gray] at (0.84,-3.58) {2};
\node[font=\tiny,gray] at (1.26,-3.58) {3};
\node[font=\tiny,gray] at (1.68,-3.58) {4};
\node[font=\tiny,gray] at (2.10,-3.58) {5};
\node[font=\tiny,gray] at (2.52,-3.58) {6};
\node[font=\tiny,gray] at (2.94,-3.58) {7};
\node[font=\tiny,gray] at (3.36,-3.58) {8};
\node[font=\tiny,gray] at (3.78,-3.58) {9};
\node[font=\tiny,gray] at (4.20,-3.58) {10};
\end{tikzpicture}&
\begin{tikzpicture}[scale=0.85]
\draw[red,thick] (0.42,-0.50)--(4.20,-0.50);
\fill[red] (0.42,-0.50) circle (1.6pt);
\fill[red] (0.84,-0.50) circle (1.6pt);
\fill[red] (1.26,-0.50) circle (1.6pt);
\fill[red] (1.68,-0.50) circle (1.6pt);
\fill[red] (2.10,-0.50) circle (1.6pt);
\fill[red] (2.52,-0.50) circle (1.6pt);
\fill[red] (2.94,-0.50) circle (1.6pt);
\fill[red] (3.36,-0.50) circle (1.6pt);
\fill[red] (3.78,-0.50) circle (1.6pt);
\fill[red] (4.20,-0.50) circle (1.6pt);
\node[left,font=\scriptsize] at (-0.25,-0.50) {$\Delta'_{5}$};
\draw[blue,thick] (3.78,-1.00)--(3.78,-1.00);
\fill[blue] (3.78,-1.00) circle (1.6pt);
\node[left,font=\scriptsize] at (-0.25,-1.00) {$2$};
\draw[blue] (2.73,-1.50) circle (1.6pt);
\node[left,font=\scriptsize] at (-0.25,-1.50) {$3$};
\draw[blue,thick] (1.68,-2.00)--(2.10,-2.00);
\fill[blue] (1.68,-2.00) circle (1.6pt);
\fill[blue] (2.10,-2.00) circle (1.6pt);
\node[left,font=\scriptsize] at (-0.25,-2.00) {$4$};
\draw[blue,thick] (1.26,-2.50)--(1.68,-2.50);
\fill[blue] (1.26,-2.50) circle (1.6pt);
\fill[blue] (1.68,-2.50) circle (1.6pt);
\node[left,font=\scriptsize] at (-0.25,-2.50) {$5$};
\draw[black,thick] (0.00,-3.00)--(0.42,-3.00);
\fill[black] (0.00,-3.00) circle (1.6pt);
\fill[black] (0.42,-3.00) circle (1.6pt);
\node[left,font=\scriptsize] at (-0.25,-3.00) {$6$};
\draw[gray!60] (-0.10,-3.40)--(4.30,-3.40);
\node[font=\tiny,gray] at (0.00,-3.58) {0};
\node[font=\tiny,gray] at (0.42,-3.58) {1};
\node[font=\tiny,gray] at (0.84,-3.58) {2};
\node[font=\tiny,gray] at (1.26,-3.58) {3};
\node[font=\tiny,gray] at (1.68,-3.58) {4};
\node[font=\tiny,gray] at (2.10,-3.58) {5};
\node[font=\tiny,gray] at (2.52,-3.58) {6};
\node[font=\tiny,gray] at (2.94,-3.58) {7};
\node[font=\tiny,gray] at (3.36,-3.58) {8};
\node[font=\tiny,gray] at (3.78,-3.58) {9};
\node[font=\tiny,gray] at (4.20,-3.58) {10};
\end{tikzpicture}&
\begin{tikzpicture}[scale=0.85]
\draw[red,thick] (0.00,-0.50)--(4.20,-0.50);
\fill[red] (0.00,-0.50) circle (1.6pt);
\fill[red] (0.42,-0.50) circle (1.6pt);
\fill[red] (0.84,-0.50) circle (1.6pt);
\fill[red] (1.26,-0.50) circle (1.6pt);
\fill[red] (1.68,-0.50) circle (1.6pt);
\fill[red] (2.10,-0.50) circle (1.6pt);
\fill[red] (2.52,-0.50) circle (1.6pt);
\fill[red] (2.94,-0.50) circle (1.6pt);
\fill[red] (3.36,-0.50) circle (1.6pt);
\fill[red] (3.78,-0.50) circle (1.6pt);
\fill[red] (4.20,-0.50) circle (1.6pt);
\node[left,font=\scriptsize] at (-0.25,-0.50) {$\Delta'_{6}$};
\draw[blue,thick] (3.78,-1.00)--(3.78,-1.00);
\fill[blue] (3.78,-1.00) circle (1.6pt);
\node[left,font=\scriptsize] at (-0.25,-1.00) {$2$};
\draw[blue] (2.73,-1.50) circle (1.6pt);
\node[left,font=\scriptsize] at (-0.25,-1.50) {$3$};
\draw[blue,thick] (1.68,-2.00)--(2.10,-2.00);
\fill[blue] (1.68,-2.00) circle (1.6pt);
\fill[blue] (2.10,-2.00) circle (1.6pt);
\node[left,font=\scriptsize] at (-0.25,-2.00) {$4$};
\draw[blue,thick] (1.26,-2.50)--(1.68,-2.50);
\fill[blue] (1.26,-2.50) circle (1.6pt);
\fill[blue] (1.68,-2.50) circle (1.6pt);
\node[left,font=\scriptsize] at (-0.25,-2.50) {$5$};
\draw[blue,thick] (0.42,-3.00)--(0.42,-3.00);
\fill[blue] (0.42,-3.00) circle (1.6pt);
\node[left,font=\scriptsize] at (-0.25,-3.00) {$6$};
\draw[gray!60] (-0.10,-3.40)--(4.30,-3.40);
\node[font=\tiny,gray] at (0.00,-3.58) {0};
\node[font=\tiny,gray] at (0.42,-3.58) {1};
\node[font=\tiny,gray] at (0.84,-3.58) {2};
\node[font=\tiny,gray] at (1.26,-3.58) {3};
\node[font=\tiny,gray] at (1.68,-3.58) {4};
\node[font=\tiny,gray] at (2.10,-3.58) {5};
\node[font=\tiny,gray] at (2.52,-3.58) {6};
\node[font=\tiny,gray] at (2.94,-3.58) {7};
\node[font=\tiny,gray] at (3.36,-3.58) {8};
\node[font=\tiny,gray] at (3.78,-3.58) {9};
\node[font=\tiny,gray] at (4.20,-3.58) {10};
\end{tikzpicture}\\
$k=4$ & $k=5$ & $k=6$
\end{tabular}
\caption{The minors of the running example: $\Delta'_k$ in red, the shifted rows in blue.}\label{fig:minors}
\end{figure}

\begin{example}\label{ex:run3} For the ladder of Example \ref{ex:run}, Figure \ref{fig:minors} shows,
for $k=1,\ldots,6$, the segment $\Delta'_k=[a_k;b_1]$ (in red) and the generalized ladder
$(\Delta^{(k)}_2,\ldots,\Delta^{(k)}_6)$ of the minor of index $k$ (the row labelled $i$ is
$\Delta^{(k)}_i$): the rows $i\leq k$ (in blue) are
the segments $[a_{i-1};b_i]$, obtained by shifting the beginnings up by one row, and the rows
$i>k$ are unchanged. No minor vanishes. For $k\geq3$, the row $i=3$ is the empty segment
$[a_2;b_3]=[7;6]$ (drawn as an open circle), so that the minors of index $3,\ldots,6$ are
generalized ladders, to which Lemma \ref{lem:empty} applies. The sets $S^0_{6,k}=\{w\in S_6^0:
w(1)=k\}$ have $24,12,4,4,2,2$ elements, and $c_k$ is their element of minimal length $k-1$.
\end{example}

\begin{lemma}\label{lem:socle} If $\pi_k\neq 0$, the socle of $\pi_k\times\Delta'_k$ is irreducible
and isomorphic to $L(\caM_{c_k})$. In particular, $\mathrm{soc}(\pi_1\times\Delta_1)=\pi$.
\end{lemma}

\begin{proof}
By \S 3.1 applied to the ladder $\pi_k$, we have $\pi_k\hookrightarrow\caJ^{(k)}_{\Id}$.
Parabolic induction is exact, so
\[\pi_k\times\Delta'_k\hookrightarrow\caJ^{(k)}_{\Id}\times\Delta'_k\simeq\caJ_{c_k}\]
by Lemma \ref{lem:minors} (iv). Hence
$\mathrm{soc}(\pi_k\times\Delta'_k)\subset\mathrm{soc}(\caJ_{c_k})=L(\caM_{c_k})$,
 and the socle of a non-zero representation of finite length is non-zero.
\end{proof}

\subsection{Words and Jacquet modules}\label{ssec:words}
Let $r_{\min}$ denote the Jacquet functor with respect to the standard parabolic subgroup of type
$(d,\ldots,d)$. For $\tau$ of finite length with cuspidal support in $\{\nu^i\rho\}$, we write
$r_{\min}(\tau)$ for its image in the Grothendieck group, viewed as a formal sum of \emph{words}
$(c_1,\ldots,c_r)$, where $(c_1,\ldots,c_r)$ stands for
$\nu^{c_1}\rho\otimes\cdots\otimes\nu^{c_r}\rho$. For two words $u=(c_1,\ldots,c_r)$ and
$u'=(c'_1,\ldots,c'_{r'})$, we denote by $u\cdot u'=(c_1,\ldots,c_r,c'_1,\ldots,c'_{r'})$ their
concatenation, extended bilinearly to formal sums of words. Let $u_1,\ldots,u_r$ be words of lengths
$n_1,\ldots,n_r$, and $n=n_1+\cdots+n_r$. For each map $f:\{1,\ldots,n\}\rightarrow\{1,\ldots,r\}$ with
$|f^{-1}(i)|=n_i$ for all $i$ (i.e.\ each position of $\{1,\ldots,n\}$ is assigned to one of the words,
the word $u_i$ receiving $n_i$ positions), the word of length $n$ whose entries at the positions of
$f^{-1}(i)$ are those of $u_i$, in the same order, is called a \emph{shuffle} of $u_1,\ldots,u_r$. The
\emph{sum of the shuffles} of $u_1,\ldots,u_r$ is the formal sum of these words over all the
$\binom{n}{n_1,\ldots,n_r}$ such maps. A given word may thus occur several times: for instance
$(1,1)$ occurs twice in the sum of the shuffles of $u_1=(1)$ and $u_2=(1)$, once for the map
with $f(1)=1$, $f(2)=2$ (the first position is assigned to $u_1$), and once for the map with $f(1)=2$,
$f(2)=1$ (it is assigned to $u_2$). We use three
standard facts.
\begin{itemize}
\item $r_{\min}(\Delta([a;b]))$ is the single word $\mathbf{w}[a;b]=(b,b-1,\ldots,a)$, empty if
$a=b+1$ (\cite{Zel80}).
\item By the geometric lemma \cite{BZ77}, $r_{\min}(\tau_1\times\cdots\times\tau_r)$ is obtained by
taking, for all words $u_1,\ldots,u_r$ of $r_{\min}(\tau_1),\ldots,r_{\min}(\tau_r)$ (with their
multiplicities), the sum of the shuffles of $u_1,\ldots,u_r$. For instance
$r_{\min}(\nu\rho\times\nu\rho)=2\,(1,1)$.
\item Jacquet functors are exact, and $r_{\min}(\tau\otimes\tau')=r_{\min}(\tau)\cdot r_{\min}(\tau')$
for an external tensor product.
\end{itemize}
For $v\in S_t^0$, the \emph{costandard word} of $\caM_v$ is the concatenation
\[W(v)=\mathbf{w}[a_{v(t)};b_t]\cdot\mathbf{w}[a_{v(t-1)};b_{t-1}]\cdot \, \cdots \, \cdot \mathbf{w}[a_{v(1)};b_1]\]
of the words of the segments of $\caM_v$, in the order of $\caJ_v$.

\begin{lemma}\label{lem:costword} For $v\in S_t^0$, the word $W(v)$ occurs in
$r_{\min}(L(\caM_v))$. Consequently, for every representation $X$ of finite length, the
multiplicity of $L(\caM_v)$ in $X$ is bounded by the multiplicity of $W(v)$ in $r_{\min}(X)$.
\end{lemma}

\begin{proof}
Since $L(\caM_v)\hookrightarrow\caJ_v$ (\S 3.1), Frobenius reciprocity shows that
$\Delta([a_{v(t)};b_t])\otimes\cdots\otimes\Delta([a_{v(1)};b_1])$ is a quotient of the
corresponding Jacquet module of $L(\caM_v)$. Applying $r_{\min}$, we get the first assertion. The
second follows from the exactness of $r_{\min}$.
\end{proof}

\begin{lemma}\label{lem:constituents} For $w\in S_t^0$, every irreducible subquotient of $\caI_w$,
or equivalently of $\caJ_w$, is of the form $L(\caM_x)$ with $x\in S_t^0$.
\end{lemma}

\begin{proof}
By \cite{Zel80}, Theorem 7.1, the irreducible subquotients of a product of Zelevinsky's
representations $\langle\Delta_1\rangle\times\cdots\times\langle\Delta_r\rangle$ are the $\langle\frn\rangle$ with
$\frn\leq(\Delta_1,\ldots,\Delta_r)$. The Zelevinsky involution is a ring automorphism of $\caR$ which
exchanges $\langle\Delta\rangle$ and $\Delta(\Delta)$ and maps $\langle\frn\rangle$ to $L(\frn)$ (\cite{Zel80}, \S9;
\cite{MW86}). Hence the irreducible subquotients of $\Delta_1\times\cdots\times\Delta_r$ are the $L(\frn)$
with $\frn\leq(\Delta_1,\ldots,\Delta_r)$ (see also \cite{LM13}, \S6.1).

Let $L(\frn)$ be an irreducible subquotient of $\caI_w$, so that $\frn\leq\caM_w$. An elementary operation
on two linked segments $[a;b]$ and $[a';b']$, with $a<a'\leq b+1\leq b'$, replaces them by $[a;b']$ and
$[a';b]$, the latter being the empty segment $[b+1;b]$ if they are juxtaposed. If we keep this empty
segment, the multisets of the beginnings and of the ends of the segments are unchanged. Starting from
$\caM_w$, with its empty segments, we obtain $t$ segments containing the non-empty segments of $\frn$,
whose beginnings are the $a_j$ and whose ends are the $b_i$, each exactly once. Writing the segment
with end $b_i$ as $[a_{x(i)};b_i]$ defines $x\in S_t$, with $a_{x(i)}\leq b_i+1$, i.e.\ $x\in S_t^0$, and
$\frn=\caM_x$.
\end{proof}

\noindent In fact $x\geq w$, since for $w,x\in S_t^0$ one has $\caM_x\leq\caM_w$ if and only if $x\geq w$
(\cite{LM13}, \S6.1); this is consistent with Proposition \ref{fact:KL}. We do not need this
refinement.

\subsection{The matching table}\label{ssec:matching}
 By Lemma \ref{lem:costword}, the multiplicity of $L(\caM_y)$ in a
representation $X$ of finite length is at most the multiplicity of the costandard word $W(y)$ in
$r_{\min}(X)$. The two main combinatorial results of this section are obtained in this way: the
multiplicity one statement of \S\ref{ssec:multone}, where $X=\caI_x$ is a product of segments, and
the key lemma (Proposition \ref{conj:key}), where $X=\pi_k\times\Delta'_k$. In both cases one has to
count how many times a given word $W(y)$ occurs in $r_{\min}(X)$, and to show that the answer is $0$
or $1$ in the relevant cases.

By the geometric lemma (\S\ref{ssec:words}), $r_{\min}$ of a product of segments is the sum of the
shuffles of their words. So we have to count the ways of producing the word $W(y)$ by shuffling the
words of the segments of the product. Such a shuffle amounts to deciding, for each entry of $W(y)$,
which segment of the product writes it; the only constraint is that each segment writes its own
entries in decreasing order, $b,b-1,\ldots,a$. The matching table records these decisions. Its cells
are the entries of $W(y)$, arranged by their value (the rows, or \emph{levels}) and by the segment of
$\caM_y$ to which they belong (the columns, or \emph{targets}); a filling indicates which segment of
the product (a \emph{factor}) writes each of them. The count is then made level by level, from the
top level down. The passage from the level $c+1$ to the level $c$ is called a \emph{transition}.
Lemma \ref{lem:chains} below shows that a transition is forced, except at the few levels where a
factor or a target ends without the other, or begins without the other, and that the choices at these
levels are described by decreasing chains of targets. In the applications only a few such levels
occur, which makes the count explicit. For $X=\pi_k\times\Delta'_k$, which is not a product of
segments, the same table is used with one more condition on the fillings, coming from the Jacquet
modules of ladders (condition (DC) in the proof of Lemma \ref{lem:words}). The terminology (factors, targets, levels, transitions,
\ldots) is ours.

Let $x,y\in S_t^0$. Write $e_m=a_{x(m)}$ and $f_k=a_{y(k)}$. The \emph{factors} are the segments
$F_m=[e_m;b_m]$ of $\caM_x$: they are the factors of the product $\caI_x$, whose words are shuffled
to produce $r_{\min}(\caI_x)$. The \emph{targets} are the segments $T_k=[f_k;b_k]$ of $\caM_y$: the
word $W(y)$ is the concatenation of their words. Note that the factor $m$ and the target $m$ have
the same top $b_m$. An integer $c$ is called a \emph{level}. Let $X_c$ (resp. $Y_c$) be the set of
indices $m$ such that $c\in F_m$ (resp. $c\in T_m$); we say that the corresponding factors (resp.
targets) are \emph{present} at the level $c$. The sets $X_c$ and $Y_c$ have the same cardinality,
since $\caM_x$ and $\caM_y$ have the same cuspidal support.

The \emph{matching table} of $(x,y)$ has one row for each level $c$, the levels decreasing from top
to bottom, and one column for each $k=1,\ldots,t$. The cell $(c,k)$ is \emph{open} if $k\in Y_c$. A
\emph{filling} of the table puts in each open cell $(c,k)$ the index of a factor present at the level
$c$, in such a way that each factor present at $c$ appears exactly once in the row $c$. Equivalently,
a filling is a family of bijections $\sigma_c:X_c\rightarrow Y_c$: the factor $m$ is put in the cell
$(c,\sigma_c(m))$, and we say that it \emph{sits in} the target $\sigma_c(m)$ at the level $c$, and
that it is the \emph{occupant} of this target. The filling is \emph{admissible} if, going down, no
factor moves to the right:
\[\sigma_c(m)\leq\sigma_{c+1}(m)\qquad\text{whenever }m\in X_c\cap X_{c+1}.\]

\begin{lemma}\label{lem:matching} The multiplicity of $W(y)$ in $r_{\min}(\caI_x)$ is the number of
admissible fillings of the matching table of $(x,y)$.
\end{lemma}

\begin{proof}
The entry $c$ occurs in $W(y)$ once in each target $k\in Y_c$. A shuffle producing $W(y)$ attaches
each such occurrence to a factor containing $c$, which defines a bijection
$\sigma_c^{-1}:Y_c\rightarrow X_c$. The targets are read in decreasing order of their index, and
each target is read in decreasing order. Hence the occurrence of $c+1$ in the target $k'$ precedes
the occurrence of $c$ in the target $k$ if and only if $k'\geq k$. The condition that each factor
reads its entries in order is therefore $\sigma_{c+1}(m)\geq\sigma_c(m)$.
\end{proof}

\begin{example}\label{ex:table} In the running example, take $x=\Id$ and $y=s_1$ (Example
\ref{ex:run2}). The factors are $[9;10]$, $[7;9]$, $[4;6]$, $[3;5]$, $[1;4]$, $[0;1]$, and the targets are
$[7;10]$, $[9;9]$, $[4;6]$, $[3;5]$, $[1;4]$, $[0;1]$. The table has a unique admissible filling:
\begin{center}
\begin{tabular}{c|cccccc}
$c\backslash k$ & $1$ & $2$ & $3$ & $4$ & $5$ & $6$\\\hline
$10$ & $1$ &&&&&\\
$9$ & $1$ & $2$ &&&&\\
$8$ & $2$ &&&&&\\
$7$ & $2$ &&&&&\\
$6$ &&& $3$ &&&\\
$5$ &&& $3$ & $4$ &&\\
$4$ &&& $3$ & $4$ & $5$ &\\
$3$ &&&& $4$ & $5$ &\\
$2$ &&&&& $5$ &\\
$1$ &&&&& $5$ & $6$\\
$0$ &&&&&& $6$
\end{tabular}
\end{center}
At the level $9$ the factor $1$ must stay in the column $1$, which forces the factor $2$ into the
column $2$; at the level $8$ the factor $1$ has exited, the target $2$ has exited, and the factor $2$
moves to the left, to the column $1$. All the other factors sit in their own targets. Hence $W(s_1)$
occurs exactly once in $r_{\min}(\caI_{\Id})$, as asserted by Proposition \ref{prop:twopos} below.
\end{example}

\subsubsection*{Reading the table downwards}
We read the family $(\sigma_c)$ from the top level down to the bottom level. The \emph{transition}
from the level $c+1$ to the level $c$ is the passage from $\sigma_{c+1}$ to $\sigma_c$; at a given
level, we count the bijections $\sigma_c$ compatible with $\sigma_{c+1}$ (i.e.\ such that no factor
present at both levels moves to the right). Going from level $c+1$ to level $c$, the following
happens.
\begin{itemize}
\item The factors with top $b_m=c$ and the targets with top $b_k=c$ \emph{enter}. A target containing
the top $b_m$ of a factor has index $\leq m$.
\item The factor with bottom $e_m=c+1$ and the target with bottom $f_k=c+1$ \emph{exit}. There is at
most one of each, since the $a_j$ are distinct.
\end{itemize}
Everything rests on the following \emph{rigidity principle}: a bijection $\tau$ of a finite set of
integers onto itself with $\tau(m)\leq m$ for all $m$ is the identity. Since factors never move to the
right, a transition at which the same targets are present before and after (up to entering factors
coming with their targets) moves nothing. The only freedom comes from the levels where a target has to
be filled from elsewhere, and this freedom is a single decreasing chain.

\begin{lemma}\label{lem:chains} Let $\sigma_{c+1}$ be given.
\begin{enumerate}
\item[(i)] If no factor and no target exits, and every entering factor $m$ comes with its target $m$,
then $\sigma_c$ is uniquely determined: $\sigma_c(m)=m$ for the entering factors and
$\sigma_c(m)=\sigma_{c+1}(m)$ otherwise.
\item[(ii)] Otherwise, there is at most one \emph{hole} $B$, a target which has to receive a factor
coming from elsewhere, and at most one \emph{source} $K$, whose factor has to go elsewhere:
\begin{itemize}
\item $B$ is the target of the exiting factor if this target does not exit, or the entering target if
it comes without its factor;
\item $K$ is the exiting target if its occupant does not exit, or the index $m$ of an entering
factor whose target $m$ is empty, viewed as a virtual target.
\end{itemize}
If $B=K$, or if there is no such pair, $\sigma_c$ is determined as in (i). If $B>K$, there is no
admissible $\sigma_c$. If $B<K$, the admissible $\sigma_c$ correspond to the decreasing chains
$K>o_1>\cdots>o_s>B$, where the $o_i$ are targets present at levels $c+1$ and $c$. Along such a
chain, the factor in $K$, or the entering factor, moves to $o_1$, the factor in $o_1$ moves to
$o_2$, …, and the factor in $o_s$ moves to $B$. All other factors stay in their targets.
\end{enumerate}
\end{lemma}

\begin{proof}
Let $D$ be the set consisting of the targets occupied at the level $c+1$ by a factor which does not
exit, and of the indices of the entering factors. For $o\in D$, let $\tau(o)$ be the target occupied
at the level $c$ by the factor sitting in $o$ at the level $c+1$, or by the entering factor $o$.
The family is admissible at this step if and only if $\tau$ is a bijection $D\rightarrow Y_c$ with
$\tau(o)\leq o$ for all $o$: for a factor present at both levels this is the monotonicity
condition, and an entering factor $m$ has top $c=b_m$, which only lies in targets of index $\leq m$.
Comparing $D$ and $Y_c$, one sees that $D\smallsetminus Y_c\subset\{K\}$ and $Y_c\smallsetminus
D\subset\{B\}$, with $D=Y_c$ in case (i) and when $B=K$. If $D=Y_c$, then $\tau$ is the identity,
by the remark above. Otherwise, the sequence $K,\tau(K),\tau(\tau(K)),\ldots$ is strictly
decreasing, by injectivity, until it reaches $B\notin D$, and $\tau$ is the identity on the
remaining elements of $D$, by the remark above. This is possible only if $B<K$. The intermediate
elements of the chain are present at both levels: an index $m\in D\cap Y_c$ not present at the
level $c+1$ is that of a factor entering with its target, so $b_m=c$, whereas $K$ is either a target
present at the level $c+1$ or the index of a factor with top $c$, so that $K<m$ in both cases,
since the $b_i$ are distinct. Hence the admissible $\tau$ are given by the decreasing chains of the statement.
\end{proof}

A factor \emph{moves}
between the levels $c+1$ and $c$ if $\sigma_c(m)\neq\sigma_{c+1}(m)$, and it is \emph{displaced} at
the level $c$ if $\sigma_c(m)\neq m$. An exit (of a factor or of a target) is \emph{clean} if the
transition at this level is as in Lemma \ref{lem:chains} (i), so that no factor moves.
A position $m$ with $x(m)=y(m)$ is called \emph{aligned}: the factor $m$ and the target $m$ are then
the same segment. In particular, as long as only aligned positions have exited, every factor sits in
its own target.

The following consequence of Lemma \ref{lem:chains} is used in all the counts below.

\begin{lemma}[Exit lemma]\label{lem:exit} Let $m$ be a factor and $k$ a target with the same bottom
$\beta$, and let $(\sigma_c)$ be an admissible filling. Then $\sigma_\beta(m)\leq k$. In particular, if
$\sigma_\gamma(m)>k$ at some level $\gamma>\beta$, then at some level $c$ with $\beta\leq c<\gamma$, a
chain of Lemma \ref{lem:chains} (ii) passes through the target occupied by the factor $m$ at the level
$c+1$.
\end{lemma}

\begin{proof}
At the transition from $\beta$ to $\beta-1$, the factor $m$ and the target $k$ exit, and they are the only
exiting factor and target, since the bottoms are distinct. If $\sigma_\beta(m)\neq k$, the hole is
$B=\sigma_\beta(m)$ and the source is $K=k$, and Lemma \ref{lem:chains} requires $B<K$. The second
assertion follows, since a factor only moves at the levels where it lies on a chain.
\end{proof}

\subsection{The key lemma}\label{ssec:key}
\begin{prop}[Key lemma]\label{conj:key} Assume $\pi_k\neq 0$. Then $\pi_k\times\Delta'_k$ has
length at most $2$, and
\[[\pi_k\times\Delta'_k]=[L(\caM_{c_k})]+[L(\caM_{c_{k+1}})]\quad\text{in } \caR,\]
where the second term is $0$ if $k=t$ or $\pi_{k+1}=0$.
\end{prop}

For $t=2$ and $k=1$, this is the classical exact sequence
 $0\rightarrow L(\Delta_1,\Delta_2)\rightarrow\Delta_2\times\Delta_1\rightarrow\Delta([a_1;b_2])\times\Delta([a_2;b_1])\rightarrow
 0$,
  and for $k=2$ it says that $\Delta([a_1;b_2])\times\Delta([a_2;b_1])$ is irreducible.

\begin{example}\label{ex:run4} In the running example, $\pi_1\times\Delta'_1=L(\Delta_2,\ldots,
\Delta_6)\times\Delta_1$ has length $2$, with constituents $\pi=L(\caM_{c_1})$ and $L(\caM_{c_2})$,
where $\caM_{c_2}=([7;10],[9;9],[4;6],[3;5],[1;4],[0;1])$ is the multisegment of Figure
\ref{fig:cover}. For $k=6$, $\pi_6\times\Delta'_6=L([9;9],[4;5],[3;4],[1;1])\times[0;10]$ is
irreducible, equal to $L(\caM_{c_6})$, since there is no minor of index $7$.
\end{example}

\begin{lemma}\label{lem:words} Assume $\pi_k\neq0$. For $x\in S_t^0$, the multiplicity of $W(x)$ in
$r_{\min}(\pi_k\times\Delta'_k)$ is $1$ if $x=c_k$, or if $x=c_{k+1}$ and $\pi_{k+1}\neq 0$, and
$0$ otherwise.
\end{lemma}

\begin{proof}
For simplicity of notation we assume that the segments $\Delta^{(k)}_i$ are non-empty. In general,
``consecutive rows'' below means consecutive among the non-empty ones, and the argument is unchanged.

\emph{Step 1: Jacquet modules of ladders.} We call the segments of a (generalized) ladder its
\emph{rows}, numbered from top to bottom as in Figures \ref{fig:ladder} and \ref{fig:minors}; two rows
are \emph{consecutive} if their indices are. Let $\sigma=L(\Delta''_1,\ldots,\Delta''_r)$ be a ladder
representation, with $\Delta''_i=\Delta([\alpha_i;\beta_i])$. By \cite{KL12}, Theorem 2.1, the
Jacquet module of $\sigma$ relative to a maximal standard parabolic subgroup is the sum of the
\[L([\gamma_1+1;\beta_1],\ldots,[\gamma_r+1;\beta_r])\otimes
L([\alpha_1;\gamma_1],\ldots,[\alpha_r;\gamma_r]),\]
over the $(\gamma_i)$ with $\alpha_i-1\leq\gamma_i\leq\beta_i$ and $\gamma_1>\cdots>\gamma_r$, of
the appropriate size, with the convention
$L(\ldots,[\alpha;\alpha-1],\ldots)=L(\ldots,\widehat{[\alpha;\alpha-1]},\ldots)$. The proof of
\cite{KL12} only uses Theorem \ref{LM} (ii), so this also holds for generalized ladders (Lemma
\ref{lem:empty}). Iterating, $r_{\min}(\sigma)$ is the sum of the words produced by the following
procedure. Start with $\gamma=\beta$. At each step, choose a row $i$ with $\gamma_i\geq\alpha_i$
and $\gamma_i-1>\gamma_{i+1}$, write the entry $\gamma_i$, and replace $\gamma_i$ by $\gamma_i-1$.
Since $\alpha_i>\alpha_{i+1}$, the condition $\gamma_i-1>\gamma_{i+1}$ means: if the row $i+1$
contains $\gamma_i-1$, it has already written it.

\emph{Step 2:} We apply this to $\pi_k$, whose rows are the segments $\Delta^{(k)}_i$,
$2\leq i\leq t$, and shuffle with the word of $\Delta'_k$. We use the matching table of
\S\ref{ssec:matching}, with $x$ replaced by $c_k$ and $y$ by $x$: the factors are the segments of
$\caM_{c_k}$ (the factor $1$ is $\Delta'_k$, and the factors $i\geq2$ are the rows of $\pi_k$), and the
targets are the segments of $\caM_x$. Since $\pi_k$ is not the product of its rows, the words of
$r_{\min}(\pi_k)$ are not all the shuffles of the words of the rows: by Step 1, two consecutive rows
must write their entries in a certain relative order. Accordingly, the multiplicity of $W(x)$ in
$r_{\min}(\pi_k\times\Delta'_k)$ is the number of fillings $(\sigma_c)$ of the table, still called
\emph{admissible}, such that:
\begin{itemize}
\item (monotonicity) $\sigma_c(m)\leq\sigma_{c+1}(m)$;
\item (DC) for consecutive rows $i<i'$, if $\ell$ is an entry of the row $i$ and $\ell-1$ an entry
of the row $i'$, then $\sigma_\ell(i)<\sigma_{\ell-1}(i')$.
\end{itemize}
The first condition is that of Lemma \ref{lem:matching}. The second one, (DC) (for ``descending
consecutive rows''), expresses the constraint of Step 1: the occurrence of $\ell-1$ read by the row
$i'$ precedes the occurrence of $\ell$ read by the row $i$; in the table, the factor $i'$ at the level
$\ell-1$ is strictly to the right of the factor $i$ at the level $\ell$.

\emph{Step 3:} The top $b_1$ of the factor $1$ only lies in the target
$1$, which has the smallest index. By monotonicity, the factor $1$ stays in the target $1$. Let
$i<i'$ be consecutive rows present at a level $c$. Then the row $i$ is present at $c+1$, since
otherwise $c=b_i$ and $b_{i'}<c$. By monotonicity and (DC),
$\sigma_c(i)\leq\sigma_{c+1}(i)<\sigma_c(i')$. Hence each $\sigma_c$ is the increasing bijection
$X_c\rightarrow Y_c$, and there is at most one admissible family, which we denote by $\theta$.

Note also that the rows present at a level $c$ form a set of consecutive rows, since the tops $b_i$
and the bottoms of the rows both decrease with $i$.

\emph{Step 4: identification of $x$.} Each value $a_j$ is the bottom of exactly one factor $s(j)$ and
of exactly one target $\beta(j)=x^{-1}(j)$; an empty factor (resp.\ target) is viewed as entering
without its target (resp.\ factor) at the level $a_j-1$. We read the levels downwards. As long as
$s(j)=\beta(j)$ for the values already passed, $\theta$ is the identity, by Lemma \ref{lem:chains} (i)
(entering factors come with their targets, since factors and targets have the same tops). At a level
$c+1=a_j$ with $s=s(j)\neq\beta=\beta(j)$, Lemma \ref{lem:chains} (ii) applies with the hole
$B=\theta_{c+1}(s)$ (or $B=s$ if the target $s$ enters without its factor) and the source $K=\beta$. If the
factor $s$ is empty, the target $s$ enters at the level $c=b_s$, whereas $\beta$ has a top $>c$, hence
$K=\beta<s=B$, which is excluded by Lemma \ref{lem:chains} (ii). Otherwise the transition is a decreasing chain from $K$ to $B$, and, since
$\theta_c$ is the increasing bijection, every factor sitting in a target of the chain moves to the
previous target of the chain, the factor $s$ leaving $B$. For consecutive rows $u<u'$, with $u$ present
at the level $c+1$ and $u'$ at the level $c$, (DC) at $\ell=c+1$ gives $\theta_{c+1}(u)<\theta_c(u')$: the
row $u'$ cannot take the target which the row $u$ occupied at the level $c+1$. The factor $1$ never moves
(Step 3), and the factors present at a level are the factor $1$, if present, and a set of consecutive
rows; so a chain with two or more moving factors, or a moving factor following a row in the hole, would
contradict (DC). We get:
\begin{itemize}
\item[$(*)$] exactly one factor moves, into the hole $B$, and it is the factor following the
occupant of $B$; if this occupant is a row, we get a contradiction.
\end{itemize}
We apply $(*)$ at the successive levels where $s(j)\neq\beta(j)$.
\begin{itemize}
\item At the first such level, the occupant of the hole is the factor $s$ in its own target. By $(*)$,
$s=1$, so $j=k$ and $c+1=a_k$. The rows $2,\ldots,k$ have bottoms $a_1,\ldots,a_{k-1}>a_k$, so the
factor following the factor $1$ is the row $k+1$: it moves from its target $\beta=k+1$ to the target
$1$, so $x(k+1)=k$, and it must contain $a_k-1$, i.e.\ $a_k\leq b_{k+1}+1$, which holds if and only if
$\pi_{k+1}\neq0$.
\item The row $k+1$ then sits in the target $1$ down to its bottom $a_{k+1}$. If the target
$\beta(k+1)$ is not the target $1$, $(*)$ at the level $a_{k+1}$ gives a contradiction, the occupant of
the hole being the row $k+1$. Hence $x(1)=k+1$.
\item At any further level with $s(j)\neq\beta(j)$, the occupant of the hole is a row, which contradicts
$(*)$.
\end{itemize}
Hence $x=c_k$, or $x=c_{k+1}$ with $\pi_{k+1}\neq0$. Conversely, in these two cases one checks
directly that $\theta$ satisfies monotonicity and (DC).
\begin{itemize}
\item For $c_k$, $\theta$ is the identity.
\item For $c_{k+1}$, $\theta$ is the identity except that the row $k+1$ sits in the target $1$ at
the levels $a_{k+1}\leq c<a_k$. There (DC) only involves the pairs $(k,k+1)$ and $(k+1,k+2)$, and
reads $k<k+1$ and $1<k+2$.
\end{itemize}
For $k=1$ and $x=c_2=s_1$ in the running example, the filling $\theta$ is the one of Example
\ref{ex:table}: the row $2$ moves to the target $1$ at the level $a_1-1=8$.
\end{proof}

\begin{proof}[Proof of Proposition \ref{conj:key}]
Let $k_0$ be the largest $k$ such that $\pi_k\neq 0$, and let $k\leq k_0$. By Lemma
\ref{lem:constituents}, the irreducible subquotients of $\pi_k\times\Delta'_k\subset\caJ_{c_k}$ are
of the form $L(\caM_x)$, $x\in S_t^0$. By Lemmas \ref{lem:costword} and \ref{lem:words}, only
$L(\caM_{c_k})$ and $L(\caM_{c_{k+1}})$ can occur, each at most once. By Lemma \ref{lem:socle},
$L(\caM_{c_k})$ occurs. Hence
\[[\pi_k\times\Delta'_k]=[L(\caM_{c_k})]+\varepsilon_k[L(\caM_{c_{k+1}})],\qquad\varepsilon_k\in\{0,1\},\]
with $\varepsilon_{k_0}=0$. On the other hand, for $1\leq j\leq k_0$, Theorem \ref{LM} (ii) applied
to $\pi_j$ (Lemma \ref{lem:empty} for the empty segments) and Lemma \ref{lem:minors} (iv) give
\begin{equation}\label{eq:minorexp}[\pi_j\times\Delta'_j]=\sum_{w\in S^0_{t,j}}(-1)^{l(w)-(j-1)}[\caI_w],\end{equation}
using that $\caR$ is commutative. Summing over $j$ and using Theorem \ref{LM} (ii) for $\pi$
(this is the expansion of Tadi\'c's determinant along its first row, see \S\ref{ssec:minors}), we get
\[[\pi]=\sum_{k=1}^{k_0}(-1)^{k-1}[\pi_k\times\Delta'_k].\]
For $k<k_0$, the coefficient of $[L(\caM_{c_{k+1}})]$ in the right-hand side is
$(-1)^{k-1}\varepsilon_k+(-1)^k$, since the $c_j$ are pairwise distinct. It must vanish, since
$c_{k+1}\neq c_1$ and $L(\caM_{c_1})=\pi$. Hence $\varepsilon_k=1$.
\end{proof}

\begin{cor}\label{cor:zigzag} Let $k_0$ be the largest $k$ such that $\pi_k\neq 0$. Then for $1\leq k\leq k_0$
 there is an exact sequence
\[0\rightarrow L(\caM_{c_k})\rightarrow\pi_k\times\Delta'_k\rightarrow L(\caM_{c_{k+1}})\rightarrow 0\]
(with $L(\caM_{c_{k_0+1}})=0$), which is non-split for $k<k_0$, and the complex
\begin{equation}\label{eq:laplace}
0\rightarrow\pi\rightarrow\pi_1\times\Delta'_1\stackrel{\partial_1}{\rightarrow}
\pi_2\times\Delta'_2\stackrel{\partial_2}{\rightarrow}\cdots\stackrel{\partial_{k_0-1}}
{\rightarrow}\pi_{k_0}\times\Delta'_{k_0}\rightarrow 0\end{equation}
is exact, where $\partial_k$ is the composition $\pi_k\times\Delta'_k\twoheadrightarrow L(\caM_{c_{k+1}})
\hookrightarrow\pi_{k+1}\times\Delta'_{k+1}$. Moreover, any non-zero morphism $\pi_k\times\Delta'_k
\rightarrow\pi_{k+1}\times\Delta'_{k+1}$ is a scalar multiple of $\partial_k$.
\end{cor}

\begin{proof}
By Lemma \ref{lem:socle} the socle is simple, so a module of length $2$ with simple socle $L(\caM_{c_k})$
 is non-split, with cosocle $L(\caM_{c_{k+1}})$. Exactness of (\ref{eq:laplace}) follows, since
 $L(\caM_{c_1})=\pi$
  by Lemma \ref{lem:socle}. For the last assertion, the image of a non-zero morphism is a quotient
  of $\pi_k\times\Delta'_k$
   and a subrepresentation of $\pi_{k+1}\times\Delta'_{k+1}$. These two modules have no common
   composition factor
    other than $L(\caM_{c_{k+1}})$, which occurs with multiplicity one in each. The image is
    therefore $L(\caM_{c_{k+1}})$,
     and the morphism factors through $\partial_k$ because $\mathrm{End}(L(\caM_{c_{k+1}}))=\bbC$.
\end{proof}

\begin{rmk}\label{rmk:tadic_ck} Put $A_k=\sum_{j\geq k}(-1)^{j-k}[\pi_j\times\Delta'_j]$. By
Proposition \ref{conj:key}, $A_k=[L(\caM_{c_k})]$. On the other hand, summing the expansions
(\ref{eq:minorexp}) of the proof of Proposition \ref{conj:key} over $j\geq k$, with the signs
$(-1)^{j-k}$, gives $A_k=\sum_{w(1)\geq k}(-1)^{l(w)-(k-1)}[\caI_w]$. Since $l(c_k)=k-1$, we get
\[[L(\caM_{c_k})]=\sum_{w\in S^0_t,\ w(1)\geq k}(-1)^{l(w)-l(c_k)}[\caI_w],\]
a determinantal formula of Tadi\'c type for $L(\caM_{c_k})$, which is not a ladder representation in general.
\end{rmk}

\subsection{Multiplicity one in intervals of length at most two}\label{ssec:multone}
The aim of this subsection is the following statement.

\begin{prop}\label{prop:multone} For $x\leq y$ in $S_t^0$ with $l(y)-l(x)\leq2$, the multiplicity
of $L(\caM_y)$ in $\caJ_x$ is $1$.
\end{prop}

It is at least $1$, since $\caJ_x\twoheadrightarrow\caJ_y\supset L(\caM_y)$ by Proposition
\ref{prop:inclusion}. By Lemma \ref{lem:costword}, and since $\caI_x$ and
$\caJ_x$ have the same image in $\caR$, it is enough to show that $W(y)$ occurs exactly once in
$r_{\min}(\caI_x)$. This is done in Propositions \ref{prop:twopos}, \ref{prop:threepos} and
\ref{prop:fourpos} below, according to the number of positions at which $x$ and $y$ differ, except for
one pattern, which is deduced by duality (Corollary \ref{cor:132}).

All these proofs follow the same pattern, in the matching table of $(x,y)$
(\S\ref{ssec:matching}). Above the first exit at a non-aligned position, every factor sits in its own
target, and the transitions are forced (Lemma \ref{lem:chains} (i)). At each exit at a non-aligned
position, Lemma \ref{lem:chains} (ii) describes the possible moves as decreasing chains of targets,
and the length condition $l(y)-l(x)\leq2$ (through the restrictions on the aligned positions listed
below) leaves at most one chain which can be completed into an admissible filling. One then checks
that this filling is admissible. Example \ref{ex:table} is the simplest instance.

\subsubsection{Equality and coverings}

\begin{prop}\label{prop:twopos} Let $u,v\in S_t^0$ with $u=v$ or $u\lessdot v$. Then $W(v)$ occurs
exactly once in $r_{\min}(\caI_u)$.
\end{prop}

\begin{proof}
We show that there is exactly one admissible family. If $u=v$, all positions are aligned, every
transition is as in Lemma \ref{lem:chains} (i), and the family is unique. Let $u\lessdot v$. By
\cite{BB05}, Lemma 2.1.4, $v=u\cdot(i,j)$ with $i<j$, $u(i)<u(j)$, and for $i<m<j$ either
$u(m)<u(i)$ or $u(m)>u(j)$. Put $\alpha=a_{u(i)}>\beta=a_{u(j)}$. The factor $i$ and the target $j$
have bottom $\alpha$, the factor $j$ and the target $i$ have bottom $\beta$. Above $\alpha$, every
factor sits in its own target. At the level $\alpha$, the factor $i$ exits from the target $i$, and
the target $j$ exits (or, if it is empty, the factor $j$ enters without its target). By Lemma
\ref{lem:chains}, the factor $j$ moves along a chain $j>o_1>\cdots>o_s>i$, whose intermediate
targets are aligned targets $m\in(i,j)$ present below $\alpha$, i.e. with $u(m)>u(j)$, hence with
bottom $<\beta$. Between $\alpha$ and $\beta$, the exits are those of aligned factors which have
not been displaced, so they are clean, and the factor $j$ stays in the target $o_1$ (or $i$ if
$s=0$). At the level $\beta$, the factor $j$ exits together with the target $i$. If $s\geq1$, the
factor $j$ is in the target $o_1>i$ at the level $\beta$, which contradicts Lemma \ref{lem:exit}. Hence $s=0$, and
the family is unique. Conversely, the family in which the factor $j$ moves directly to the target
$i$ at the level $\alpha$, all other transitions being as in Lemma \ref{lem:chains} (i), is admissible.
\end{proof}

\subsubsection{Three positions} Assume that $l(y)-l(x)=2$ and that $x$ and $y$ differ exactly at
three positions $p<q<r$. Let $v_1<v_2<v_3$ be the values of $x$ on $\{p,q,r\}$, and put
$A_i=a_{v_i}$, so that $A_1>A_2>A_3$. The restrictions of $x$ and $y$ to $\{p,q,r\}$ have one of the patterns
\[123\rightarrow312,\qquad 123\rightarrow231,\qquad 132\rightarrow321,\qquad 213\rightarrow321,\]
the only intervals of length $2$ of $S_3$ whose ends differ at all three positions. A direct count
of inversions shows that each aligned position $m$ contributes $0$ or $2$ to $l(y)-l(x)-2$. Hence
all contributions vanish. Concretely, aligned positions $m$ outside $[p,r]$ are unrestricted, while
for $p<m<q$ or $q<m<r$ the value $x(m)$ must lie in the following intervals:
\begin{center}
\begin{tabular}{c|c|c}
 & $p<m<q$ & $q<m<r$\\ \hline
$123\rightarrow312$ & $<v_1$, $>v_3$ & $<v_1$, $(v_1,v_2)$, $>v_3$\\
$123\rightarrow231$ & $<v_1$, $(v_2,v_3)$, $>v_3$ & $<v_1$, $>v_3$\\
$132\rightarrow321$ & $<v_1$, $>v_3$ & $<v_1$, $(v_2,v_3)$, $>v_3$\\
$213\rightarrow321$ & $<v_1$, $(v_1,v_2)$, $>v_3$ & $<v_1$, $>v_3$
\end{tabular}
\end{center}
Recall that $x(m)<v_1$ means $e_m>A_1$, and $x(m)>v_3$ means $e_m<A_3$. Only three non-aligned
exits occur, at the levels $A_1>A_2>A_3$. The aligned exits are clean, i.e. covered by Lemma
\ref{lem:chains} (i), as long as the corresponding factor has not been displaced.

\begin{prop}\label{prop:threepos} In this situation, except for the pattern $132\rightarrow321$,
the word $W(y)$ occurs exactly once in $r_{\min}(\caI_x)$.
\end{prop}

\begin{proof}
It occurs at least once, since $L(\caM_y)$ is a subquotient of $\caJ_x$. We show that there is at
most one admissible family. Above $A_1$ every factor sits in its own target.

\emph{Pattern $123\rightarrow312$.} The factors $p,q,r$ have bottoms $A_1,A_2,A_3$, and the targets
$p,q,r$ have bottoms $A_3,A_1,A_2$.
\begin{itemize}
\item At $A_1$, the factor $p$ exits from its target $p$, and the target $q$ exits; if the target $q$
is empty, the factor $q$ enters instead. By Lemma \ref{lem:chains}, the factor $q$ moves down along
a chain from $q$ to $p$. The possible intermediate targets are the aligned $m\in(p,q)$ present below
$A_1$, i.e. with $x(m)>v_3$.
\item Between $A_1$ and $A_2$, only aligned factors which have not been displaced exit.
\item At $A_2$, the factor $q$ exits, and the target $r$ exits (or the factor $r$ enters). The
factor $r$ moves down along a chain ending at the target of the factor $q$, which is $\geq p$.
\item Between $A_2$ and $A_3$, no aligned exit occurs in $(p,r)$, since no aligned value lies in
$(v_2,v_3)$ there, and the exits outside $[p,r]$ do not move anything.
\item At $A_3$, the factor $r$ exits together with the target $p$, so it sits in the target $p$
by Lemma \ref{lem:exit}, its target being $\geq p$. This happens only if both chains above have no
intermediate target.
\end{itemize}
Hence there is a unique admissible family.

\emph{Pattern $123\rightarrow231$.} The factors $p,q,r$ have bottoms $A_1,A_2,A_3$, and the targets
$p,q,r$ have bottoms $A_2,A_3,A_1$.
\begin{itemize}
\item At $A_1$, the factor $p$ exits from the target $p$, and the target $r$ exits. Since $A_1\leq
b_r+1\leq b_q$, the target $q$ is present below $A_1$ and is a possible intermediate target, besides
the aligned targets of $(p,q)$ and $(q,r)$ with values $>v_1$.
\item At $A_2$, the target $p$ exits together with the factor $q$. Since no target of index $<p$ is
involved in the chains, the factor $q$ must then be in the target $p$, by Lemma \ref{lem:exit}. Between $A_1$ and $A_2$ no
aligned position of $(p,r)$ exits, since no aligned value lies in $(v_1,v_2)$ there. Hence the
chain at $A_1$ must contain $q$, and its part below $q$ must be empty, so that the factor $q$ moves
directly to $p$.
\item If the chain at $A_1$ contains an intermediate target $m\in(q,r)$, the factor $r$ ends in a
target $>q$. Between $A_2$ and $A_3$ the only aligned exits in $(p,r)$ have values in $(v_2,v_3)$,
hence lie in $(p,q)$; they are clean and do not move the factor $r$. At $A_3$, the factor $r$ exits
together with the target $q$, while sitting in a target $>q$, which contradicts Lemma \ref{lem:exit}.
\end{itemize}
Hence the chain at $A_1$ is $r>q>p$, and the family is unique.

\emph{Pattern $213\rightarrow321$.} The factors $p,q,r$ have bottoms $A_2,A_1,A_3$, and the targets
$p,q,r$ have bottoms $A_3,A_2,A_1$.
\begin{itemize}
\item At $A_1$, the factor $q$ exits from its target, and the target $r$ exits (or the factor $r$
enters). The factor $r$ moves down to a target of $(q,r)$, or to $q$ if the chain has no
intermediate target. The intermediate targets are aligned targets of $(q,r)$.
\item Between $A_1$ and $A_2$, the aligned exits in $(p,r)$ have values in $(v_1,v_2)$, hence lie
in $(p,q)$, and are clean.
\item At $A_2$, the factor $p$ exits from its target $p$, and the target $q$ exits. The chain goes
from $q$ down to $p$ through aligned targets of $(p,q)$. It does not move a factor sitting in a
target of $(q,r)$.
\item Between $A_2$ and $A_3$ there is no aligned exit in $(p,r)$.
\item At $A_3$, the factor $r$ exits together with the target $p$. By Lemma \ref{lem:exit} the factor
$r$ is then in the target $p$, i.e.\ both chains have no intermediate target.
\end{itemize}
Hence the family is unique.

\end{proof}

For the pattern $132\rightarrow321$, the chain at $A_1$ may pass through the target $q$ or not, and
$W(y)$ may occur more than once. We use duality instead.

\begin{cor}\label{cor:132} In the situation of Proposition \ref{prop:threepos}, with the pattern
$132\rightarrow321$, the multiplicity of $L(\caM_y)$ in $\caJ_x$ is $1$.
\end{cor}

\begin{proof} Consider the dual ladder $\pi^\vee$ and the elements $\tilde x=w_0x^{-1}w_0$, $\tilde
y=w_0y^{-1}w_0$ (see the introduction). Taking contragredients, the multiplicity of $L(\caM_y)$ in
$\caJ_x$ is that of $L(\caM_y)^\vee=L(\caM'_{\tilde y})$ in $\caJ_x^\vee\simeq\caI'_{\tilde x}$. One
checks that $\tilde x$ and $\tilde y$ differ at three positions, with the pattern $213\rightarrow321$.
By Proposition \ref{prop:threepos} and Lemma \ref{lem:costword}, this multiplicity is at most $1$, and
it is at least $1$ as above.
\end{proof}
\subsubsection{Four positions} Assume now that $l(y)-l(x)=2$ and that $x$ and $y$ differ exactly at
four positions $p_1<p_2<p_3<p_4$. Let $v_1<\cdots<v_4$ be the values of $x$ on these positions. The
restrictions of $x$ and $y$ to $\{p_1,\ldots,p_4\}$ have one of the eleven patterns of the
following table, which are the intervals of length $2$ of $S_4$ whose ends differ at all four
positions. As in the case of three positions, a count of inversions shows that each aligned
position $m$ contributes a non-negative even number to $l(y)-l(x)-2$, hence $0$. This puts no
restriction on the aligned positions outside $[p_1,p_4]$, and for $p_s<m<p_{s+1}$ it means that
$x(m)$ does not lie in the intervals listed in the column $(p_s,p_{s+1})$. The second column gives
the decomposition $y=x\cdot(i_1,j_1)(i_2,j_2)$ into two disjoint pairs $P_s=(i_s,j_s)$, where the
positions are numbered $1,\ldots,4$ in the order of $p_1,\ldots,p_4$. In all cases $x(i_s)<x(j_s)$.
\begin{center}
{\footnotesize\setlength{\tabcolsep}{3pt}
\begin{tabular}{c|c|c|c|c}
pattern & pairs & $(p_1,p_2)$ & $(p_2,p_3)$ & $(p_3,p_4)$\\ \hline
$1234\rightarrow2143$ & $(1,2)\,(3,4)$ & $(v_1,v_2)$ & -- & $(v_3,v_4)$\\
$1324\rightarrow2413$ & $(1,3)\,(2,4)$ & $(v_1,v_2)$ & $(v_1,v_2), (v_3,v_4)$ & $(v_3,v_4)$\\
$1324\rightarrow3142$ & $(1,2)\,(3,4)$ & $(v_1,v_2), (v_2,v_3)$ & -- & $(v_2,v_3), (v_3,v_4)$\\
$1342\rightarrow2431$ & $(1,4)\,(2,3)$ & $(v_1,v_2)$ & $(v_1,v_2), (v_3,v_4)$ & $(v_1,v_2)$\\
$1423\rightarrow4132$ & $(1,2)\,(3,4)$ & $(v_1,v_2), (v_2,v_3), (v_3,v_4)$ & -- & $(v_2,v_3)$\\
$2143\rightarrow3412$ & $(1,4)\,(2,3)$ & $(v_2,v_3)$ & $(v_1,v_2), (v_2,v_3), (v_3,v_4)$ & $(v_2,v_3)$\\
$2314\rightarrow3241$ & $(1,2)\,(3,4)$ & $(v_2,v_3)$ & -- & $(v_1,v_2), (v_2,v_3), (v_3,v_4)$\\
$2413\rightarrow4231$ & $(1,2)\,(3,4)$ & $(v_2,v_3), (v_3,v_4)$ & -- & $(v_1,v_2), (v_2,v_3)$\\
$3124\rightarrow4213$ & $(1,4)\,(2,3)$ & $(v_3,v_4)$ & $(v_1,v_2), (v_3,v_4)$ & $(v_3,v_4)$\\
$3142\rightarrow4231$ & $(1,3)\,(2,4)$ & $(v_3,v_4)$ & $(v_1,v_2), (v_3,v_4)$ & $(v_1,v_2)$\\
$3412\rightarrow4321$ & $(1,2)\,(3,4)$ & $(v_3,v_4)$ & -- & $(v_1,v_2)$
\end{tabular}}
\end{center}
By inspection of the table, we get:
\begin{itemize}
\item[(AI)] if $m$ is aligned and $i_s<m<j_s$, then $x(m)\notin(x(i_s),x(j_s))$;
\item[(Sep)] for each pair $P=P_s$, with $P'$ the other pair, one of the following holds:
\begin{itemize}
\item[(I)] no value of $x$ on $P'$ lies in $(x(i_s),x(j_s))$;
\item[(II)] both positions of $P'$ lie on the same side of $[i_s,j_s]$;
\item[(III)] the local pattern is $2143\rightarrow3412$ and $P$ is the inner pair.
\end{itemize}
\end{itemize}
For a pair $P=(i,j)$, put $\alpha_P=a_{x(i)}>\beta_P=a_{x(j)}$. The factor $i$ and the target $j$
have bottom $\alpha_P$, and the factor $j$ and the target $i$ have bottom $\beta_P$. By (AI), an
aligned target $m$ with $i<m<j$ which is present below $\alpha_P$ has bottom $<\beta_P$.

\begin{prop}\label{prop:fourpos} In this situation, $W(y)$ occurs exactly once in $r_{\min}(\caI_x)$.
\end{prop}

\begin{proof}
Let $H_0$ be the following family. At the level $\alpha_P$ of each pair, the chain is $j>i$: the
factor $j$ moves to the target $i$. All other transitions are as in Lemma \ref{lem:chains} (i). It
is admissible, since at $\beta_P$ the factor $j$ sits in the target $i$, which exits with it. Let
$H$ be another admissible family, and consider the highest level at which $H$ differs from $H_0$.
Up to that level, the factors of $H$ sit where they sit in $H_0$. At this point all exits other
than the $\alpha_P$ are clean, so the first difference is a chain $j>o_1>\cdots>o_s>i$ with
$s\geq1$ at some level $\alpha_P$, $P=(i,j)$. The factor $j$ is then in the target $o_1\in(i,j)$. Since the
factor $j$ and the target $i$ have the same bottom $\beta_P$, Lemma \ref{lem:exit} requires a chain, at a
level of the window $[\beta_P,\alpha_P)$, passing through the target of the factor $j$. We show that
there is none.
\begin{itemize}
\item The aligned factors displaced at $\alpha_P$ have bottoms $<\beta_P$ by (AI), and the
non-aligned ones exit at levels of the other pair $P'$.
\item \emph{Case (I).} No exit of $P'$ lies in the window. Since $H$ and $H_0$ agree above
$\alpha_P$, a chain of $P'$ at a level above $\alpha_P$ has no intermediate target and displaces no
aligned factor, and the aligned factors displaced at $\alpha_P$ exit below $\beta_P$ by (AI). Hence
the aligned exits in the window concern factors which have not been displaced, so they are clean,
and no chain occurs in the window.
\item \emph{Case (II).} The exits of $P'=(i',j')$ in the window, and the exits of aligned factors
displaced by chains of $P'$, give chains contained in $[i',j']$. Indeed, the factors of $P'$ and
the factors displaced by $P'$ only move inside $[i',j']$. This interval is disjoint from $(i,j)$,
so these chains do not move the factor $j$. The other aligned exits in the window are clean.
\item \emph{Case (III).} The local pattern is $2143\rightarrow3412$ and $P=(i,j)$ is the inner
pair, so $\alpha_P=A_1$ and $\beta_P=A_4$, with the notation $A_1>\cdots>A_4$ of the values. The
outer pair $P'=(i',j')$ has $\alpha_{P'}=A_2$, $\beta_{P'}=A_3$, and $i'<i<j<j'$.
\begin{itemize}
\item At $A_2$, the chain goes from $j'$ to $i'$.
\item At $A_3$, the factor $j'$ exits together with the target $i'$. No aligned exit of $(i',j')$
lies in $(A_3,A_2)$ for this pattern, so the factor $j'$ cannot move in between.
\item Hence the chain at $A_2$ must be $j'>i'$, without intermediate target, and it does not move
the factor $j$. The aligned exits of $(A_4,A_3)$ in $(i',j')$ are clean, since the aligned factors
displaced at $A_1$ have bottoms $<A_4$.
\end{itemize}
\end{itemize}
In all cases the factor $j$ is still in a target $>i$ at the level $\beta_P$, which contradicts Lemma
\ref{lem:exit}.
Hence $H=H_0$.
\end{proof}

\begin{proof}[Proof of Proposition \ref{prop:multone}]
If $x=y$ or $x\lessdot y$, apply Proposition \ref{prop:twopos}. If $l(y)-l(x)=2$, then $y=x\cdot
t_1t_2$ for two transpositions $t_1,t_2$, by the chain property of the Bruhat order, so $x$ and $y$
differ at most at four positions. They do not differ at exactly two positions, since a
transposition changes the length by an odd number. Hence they differ at three or four positions.
Apply Proposition \ref{prop:threepos} (with Lemma \ref{lem:costword}), Corollary \ref{cor:132} or
Proposition \ref{prop:fourpos}.
\end{proof}

\section{Construction of the coresolution and proof of Theorem \ref{BGG}}\label{sec:complex}
In this section we deduce Theorem \ref{BGG} from the results of \S\ref{sec:reduction} by homological
algebra. The only inputs are Proposition \ref{prop:inclusion} (surjections between costandard
modules), Lemma \ref{lem:minors} (combinatorics of the minors), the key lemma (Proposition
\ref{conj:key} and Corollary \ref{cor:zigzag}) and the multiplicity one statement (Proposition
\ref{prop:multone}).

\subsection{Construction of the complex}\label{ssec:complex}
Let $\pi$ be a generalized ladder representation. By the end of the proof of Proposition
\ref{prop:inclusion}, for every covering
 $w\lessdot w'$ in $S_t^0$ there is a surjection $\caJ_w\twoheadrightarrow\caJ_{w'}$.

\begin{lemma}\label{lem:rigid} Let $w\leq w'$ in $S_t^0$ with $l(w')-l(w)\leq 2$.
\begin{enumerate}
\item[(i)] The multiplicity of $L(\caM_{w'})$ in $\caJ_w$ is $1$.
\item[(ii)] $\mathrm{End}(\caJ_w)=\bbC$.
\item[(iii)] Any two surjections $\caJ_w\twoheadrightarrow\caJ_{w'}$ have the same kernel, and they
differ by a non-zero scalar.
\end{enumerate}
\end{lemma}

\begin{proof}
(i) is Proposition \ref{prop:multone}.

(ii) The socle of $\caJ_w$ is $L(\caM_w)$, which occurs once in $\caJ_w$. An endomorphism $\alpha$ acts
 on it by a scalar $\lambda$. If $\beta=\alpha-\lambda$ were non-zero, its image would be a
 non-zero submodule of $\caJ_w$, hence would contain the socle $L(\caM_w)$. Since $\beta$ kills the
 socle, $L(\caM_w)$ would occur both in $\ker\beta$ and in
 $\caJ_w/\ker\beta\simeq\mathrm{im}\,\beta$, hence at least twice in $\caJ_w$, which is absurd.
 Hence $\alpha=\lambda$.

(iii) Let $f,g:\caJ_w\twoheadrightarrow\caJ_{w'}$ be surjections. The image of $\ker g$ in
$\caJ_w/\ker f\simeq\caJ_{w'}$
 is either $0$ or contains the socle $L(\caM_{w'})$. In the second case $L(\caM_{w'})$ would occur both in
 $\ker g$ and in $\caJ_w/\ker g\simeq\caJ_{w'}$, contradicting (i). Hence $\ker g\subset\ker f$,
 and $\ker f=\ker g$
 by symmetry. Then $g=\alpha\circ f$ for an automorphism $\alpha$ of $\caJ_{w'}$, which is a scalar by (ii).
\end{proof}

For $w\in S_t^0$, choose a maximal chain $\Id=x_0\lessdot x_1\lessdot\cdots\lessdot x_m=w$ in $S_t^0$
and surjections $\caJ_{x_{r-1}}\twoheadrightarrow\caJ_{x_r}$. Let $\caK_w\subset\caJ_{\Id}$ be the kernel
 of the composite surjection $\caJ_{\Id}\twoheadrightarrow\caJ_w$.

\begin{lemma}\label{lem:kernel} The submodule $\caK_w$ does not depend on the choices. If $w\lessdot w'$,
 then $\caK_w\subset\caK_{w'}$.
\end{lemma}

\begin{proof}
By Lemma \ref{lem:rigid} (iii), the surjections along a fixed chain are unique up to scalars.
Two maximal chains of $[\Id,w]$ can be joined by a sequence of maximal chains in which two consecutive
 chains differ in exactly one element. Indeed, the order complex of $[\Id,w]$ is shellable
 (\cite{BB05}, Theorem 2.7.5),
 and in a shelling each facet meets a previous one in a face of codimension one (\cite{BB05},
 Definition A2.4.1).
 When two chains differ only in one element, the two different elements $a,b$ lie in a length-two
 interval $[x,y]$ of the chains.
  By Lemma \ref{lem:rigid} (iii), the composites
  $\caJ_x\twoheadrightarrow\caJ_a\twoheadrightarrow\caJ_y$ and
  $\caJ_x\twoheadrightarrow\caJ_b\twoheadrightarrow\caJ_y$ differ by a scalar, so the total kernel
  does not change.
  The second assertion follows by choosing a chain through $w$.
\end{proof}

By Lemma \ref{lem:kernel}, we may identify $\caJ_w$ with $\caJ_{\Id}/\caK_w$. For $w\lessdot w'$, let
$\phi_{w,w'}:\caJ_w\twoheadrightarrow\caJ_{w'}$
be the canonical projection. For every length-two interval $[x,y]=\{x,a,b,y\}$ of $S_t^0$
(\cite{BB05}, Lemma 2.7.3),
we have $\phi_{a,y}\circ\phi_{x,a}=\phi_{b,y}\circ\phi_{x,b}$.

\begin{defi} A \emph{balanced family} on $S_t^0$ is a map $s$ from the set of coverings $w\lessdot
w'$ in $S_t^0$
to $\bbC^\times$ such that $s(x,a)s(a,y)+s(x,b)s(b,y)=0$ for every length-two interval
$[x,y]=\{x,a,b,y\}$ in $S_t^0$.
For such a family, let $C^\bullet(\pi,s)$ be the complex with
\[C^j(\pi,s)=\bigoplus_{w\in S_t^0,\ l(w)=j}\caJ_w,\qquad d=\sum_{w\lessdot w'}s(w,w')\,\phi_{w,w'},\]
augmented by $\pi=\mathrm{soc}(\caJ_{\Id})\hookrightarrow C^0(\pi,s)$.
\end{defi}

By construction, and since $S_t^0$ is an interval, $d\circ d=0$. The augmentation lands in the
kernel of $d$: for $\Id\lessdot w$, if the image of $\pi$ in $\caJ_w$ were non-zero, it would be
isomorphic to $\pi$, hence to the socle $L(\caM_w)$ of $\caJ_w$, and $\caM_w=\caM$, i.e. $w=\Id$,
by Proposition \ref{prop:zelorder}. Balanced families with values in $\{\pm1\}$
exist by \cite{BB05}, Corollary 2.7.14, applied to the interval $S_t^0=[\Id,w'_0]$.

\subsection{End of the proof}

\begin{prop}\label{prop:H} Let $s$ be a balanced family on $S_t^0$, and set
$F^kC^\bullet=\bigoplus_{w(1)\geq k}\caJ_w\subset C^\bullet(\pi,s)$.
\begin{enumerate}
\item[(H1)] $(F^kC^\bullet)_k$ is a finite decreasing filtration of $C^\bullet(\pi,s)$ by subcomplexes.
\item[(H2)] If $\pi_k\neq 0$, there is a balanced family $s_k$ on the interval attached to
$(\Delta^{(k)}_2,\ldots,\Delta^{(k)}_t)$ and an isomorphism of complexes
\[F^kC^\bullet/F^{k+1}C^\bullet\simeq\big(C^\bullet(\pi_k,s_k)\times\Delta'_k\big)[-(k-1)].\]
If $\pi_k=0$, this quotient is $0$.
\item[(H3)] Assume $\pi_{k+1}\neq 0$. The component
$s(c_k,c_{k+1})\phi_{c_k,c_{k+1}}:\caJ_{c_k}\rightarrow\caJ_{c_{k+1}}$
of $d$ does not vanish on the submodule
$H^{0}(C^\bullet(\pi_k,s_k))\times\Delta'_k\simeq\pi_k\times\Delta'_k$ of $\caJ_{c_k}$.
\end{enumerate}
\end{prop}

\begin{proof}
(H1) follows from Lemma \ref{lem:minors} (v).

(H2) By Lemma \ref{lem:minors} (v), $S^0_{t,k}$ is convex in $S^0_t$. By Lemma \ref{lem:minors}
(iv), $w\mapsto\bar w$
 is an isomorphism of posets onto the interval attached to the minor, shifting lengths by $k-1$.
 Hence it identifies coverings
  and length-two intervals. For a covering $w\lessdot w'$ in $S^0_{t,k}$, both $\phi_{w,w'}$ and
   $\phi^{(k)}_{\bar w,\bar w'}\times\mathrm{id}_{\Delta'_k}$ are surjections
   $\caJ_w\twoheadrightarrow\caJ_{w'}$,
   under the identifications $\caJ_w\simeq\caJ^{(k)}_{\bar w}\times\Delta'_k$. By Lemma \ref{lem:rigid} (iii),
    they differ by a scalar $\lambda(w,w')\neq 0$. The differential of $F^k/F^{k+1}$ is then
     $\sum s(w,w')\lambda(w,w')\,\phi^{(k)}_{\bar w,\bar w'}\times\mathrm{id}$. Since its square vanishes
      and the maps $\phi^{(k)}$ commute on length-two intervals, the family $s_k(\bar w,\bar
      w')=s(w,w')\lambda(w,w')$ is balanced.

(H3) We have $c_k\lessdot c_{k+1}$ by Lemma \ref{lem:minors} (vi). By Lemma \ref{lem:rigid} (i)
 applied to $(c_k,c_{k+1})$ and to $(c_{k+1},c_{k+1})$, the kernel of $\phi_{c_k,c_{k+1}}$ does not contain
 $L(\caM_{c_{k+1}})$ as a composition factor. On the other hand, $L(\caM_{c_{k+1}})$ is a composition factor
 of $\pi_k\times\Delta'_k$ by Proposition \ref{conj:key}. Hence
 $\pi_k\times\Delta'_k\not\subset\ker\phi_{c_k,c_{k+1}}$.
\end{proof}

\begin{thm}\label{thm:main} Let $\pi$ be a generalized ladder representation and $s$ a balanced
family on $S_t^0$.
 Then $H^j(C^\bullet(\pi,s))=0$ for $j>0$, and the augmentation identifies $\pi$ with
 $H^0(C^\bullet(\pi,s))$. In particular,
 Theorem \ref{BGG} holds.
\end{thm}

\begin{proof}
By Lemma \ref{lem:empty}, we may argue by induction on the number $t$ of segments of the generalized ladder.
The case $t=1$ is clear. Consider the spectral sequence of the filtration (H1) of Proposition \ref{prop:H}. By (H2), the
induction hypothesis,
and the exactness of parabolic induction, $H^\bullet(F^k/F^{k+1})$ is $\pi_k\times\Delta'_k$, placed in degree
 $k-1$, if $\pi_k\neq0$, and $0$ otherwise. It is realized as a submodule of $\caJ_{c_k}$, the unique term of
 $F^k/F^{k+1}$ in degree $k-1$ (Lemma \ref{lem:minors} (vi)). The $E_1$ term is therefore
 concentrated on one line, and it is the complex
\[\pi_1\times\Delta'_1\stackrel{d_1}{\rightarrow}\pi_2\times\Delta'_2\stackrel{d_1}{\rightarrow}\cdots\rightarrow\pi_{k_0}\times\Delta'_{k_0}.\]
Since $c_k$ (resp. $c_{k+1}$) is the only element of length $k-1$ (resp. $k$) with $w(1)\geq k$
(resp. $w(1)\geq k+1$) by Lemma \ref{lem:minors} (vi),
 the map $d_1$ from the $k$-th term to the $(k+1)$-th term is the restriction of
 $s(c_k,c_{k+1})\phi_{c_k,c_{k+1}}$.
  It is non-zero by (H3), so it is a non-zero multiple of $\partial_k$ by Corollary \ref{cor:zigzag}.
  By the exactness of (\ref{eq:laplace}), $E_2$ is $\pi$ in degree $0$ and vanishes elsewhere.
  Hence the spectral sequence degenerates, $H^j(C^\bullet)=0$ for $j>0$, and $H^0(C^\bullet)$ has
  the same image as $\pi$ in $\caR$. Finally, $H^0(C^\bullet)\subset\caJ_{\Id}$ contains the image of
  the augmentation, which is $\pi$. Hence $H^0(C^\bullet)=\pi$.
\end{proof}

\begin{example}\label{ex:run5} In the running example, the filtration (H1) has six graded pieces.
The following table gives, for each $k$, the number of elements $w\in S^0_{6,k}$ of each length,
i.e.\ the number of costandard modules of each degree in $F^kC^\bullet/F^{k+1}C^\bullet$. By
(H2), the row $k$ is the complex of the minor $\pi_k$, shifted by $k-1$.
\begin{center}
\begin{tabular}{c|ccccccc|c}
$k\backslash j$ & $0$ & $1$ & $2$ & $3$ & $4$ & $5$ & $6$ & $|S^0_{6,k}|$\\\hline
$1$ & $1$ & $4$ & $7$ & $7$ & $4$ & $1$ & & $24$\\
$2$ & & $1$ & $3$ & $4$ & $3$ & $1$ & & $12$\\
$3$ & & & $1$ & $2$ & $1$ & & & $4$\\
$4$ & & & & $1$ & $2$ & $1$ & & $4$\\
$5$ & & & & & $1$ & $1$ & & $2$\\
$6$ & & & & & & $1$ & $1$ & $2$\\\hline
total & $1$ & $5$ & $11$ & $14$ & $11$ & $5$ & $1$ & $48$
\end{tabular}
\end{center}
The $E_1$ term of the spectral sequence is the complex
$\pi_1\times\Delta'_1\rightarrow\pi_2\times\Delta'_2\rightarrow\cdots\rightarrow\pi_6\times\Delta'_6$,
and the augmented sequence
\[0\rightarrow\pi\rightarrow\pi_1\times\Delta'_1\rightarrow\pi_2\times\Delta'_2\rightarrow\cdots
\rightarrow\pi_6\times\Delta'_6\rightarrow0\]
is exact. In this sequence, by Proposition \ref{conj:key}, the $k$-th term has the constituents $L(\caM_{c_k})$ and
$L(\caM_{c_{k+1}})$ ($k<6$), and the last one is $L(\caM_{c_6})$.
\end{example}

\section{Kazhdan-Lusztig polynomials and Cousin resolutions}\label{sec:KL}

\subsection{Relation with Kazhdan-Lusztig polynomials}\label{ssec:KL}
The multiplicities of the irreducible representations $L(\caM_v)$ in the standard modules are given
by Kazhdan-Lusztig polynomials. This gives another proof of Propositions \ref{conj:key} and
\ref{prop:multone}, which we sketch.

\begin{prop}[Zelevinsky, Chriss-Ginzburg, Arakawa-Suzuki]\label{fact:KL} For $u,v\in S_t^0$,
 the multiplicity of $L(\caM_v)$ in $\caI_u$ is $P_{u,v}(1)$, where $P_{u,v}$ is the Kazhdan-Lusztig
 polynomial of $S_t$. In particular it vanishes unless $u\leq v$.
\end{prop}

\begin{proof}[Outline of proof]
This is a special case of the $p$-adic analogue of the Kazhdan-Lusztig conjecture, conjectured by
Zelevinsky \cite{Zel81, Zel85} and proved in \cite{CG97}, Theorem 8.6.23. In the regular situation
 of ladders, the most direct route is the Arakawa-Suzuki functor, as in \cite{LM18}, \S 10,
 and \cite{BC15}, \S 4. We outline it, omitting the details of Step 1, which are as in \cite{BC15},
\S3 (types for the inertial class of $\rho^{\otimes\ell}$, compatibility with parabolic induction,
and matching of the standard modules and of the Langlands quotients). The normalization of Step 3 is
forced by Tadi\'c's formula: since $P_{e,3412}=1+q$, the other possible normalization
$P_{w_0v,w_0u}$ would contradict it for the Speh representations with $t\geq4$.

\emph{Step 1 (reduction to the graded Hecke algebra).} Let $\ell=\sum_i(b_i-a_i+1)$, so that $\pi$
and all the $\caI_w$ are representations of $G_{\ell d}$ with cuspidal support in $\{\nu^i\rho\}$.
By the theory of types \cite{BK93}, the corresponding Bernstein block is equivalent to the category
of modules over the affine Hecke algebra of $GL_\ell$. This equivalence is compatible with
parabolic induction.
By Lusztig's reduction \cite{Lus89}, the finite length modules with the relevant (real, integral)
central character
are in turn identified with finite-dimensional modules over the graded Hecke algebra $\mathbb{H}_\ell$ of
 $GL_\ell$. Under these equivalences, $\Delta([a;b])$ corresponds to the Steinberg-type module
 attached to $[a;b]$,
 the module $\caI_u$ corresponds to the standard module $\mathcal{M}(\lambda,\mu)$ of \cite{BC15}, (4.2.4),
  with the same segments, and $L(\caM_u)$ corresponds to its unique simple quotient
  $\mathcal{L}(\lambda,\mu)$.

\emph{Step 2 (Arakawa-Suzuki functor).} Let $\mathfrak{g}=\mathfrak{gl}_t$, let $\rho$ also denote
its half-sum
 of positive roots, and put
\[\lambda+\rho=(b_1+1,\ldots,b_t+1),\qquad \mu_u+\rho=(a_{u(1)},\ldots,a_{u(t)})\quad(u\in S_t).\]
By \cite{BC15}, (4.2.3), the segments of $\mathcal{M}(\lambda,\mu_u)$ are then $[a_{u(i)};b_i]$,
i.e. those of $\caM_u$.
Since $b_1>\cdots>b_t$, the weight $\lambda+\rho$ is dominant and regular. Hence the Arakawa-Suzuki functor
 $F=F_{\lambda,\mathcal{V}_\ell}$, $\mathcal{V}_\ell=(\mathbb{C}^t)^{\otimes\ell}$, from the category $\caO$
 of $\mathfrak{g}$ to $\mathbb{H}_\ell$-modules is exact (\cite{BC15}, Remark 4.1.1, \cite{AS98}).
 By \cite{BC15}, Theorem 4.2.3, condition (4.2.5) being empty because $\lambda+\rho$ is regular,
 we have $F(M(\mu_u))=\mathcal{M}(\lambda,\mu_u)$ and $F(L(\mu_u))=\mathcal{L}(\lambda,\mu_u)$
 for $u\in S_t^0$. For $u\notin S_t^0$, $\lambda-\mu_u$ is not a weight of $\mathcal{V}_\ell$, and
  both sides vanish. Therefore, for $u,v\in S_t^0$,
\[[\caI_u:L(\caM_v)]=[\mathcal{M}(\lambda,\mu_u):\mathcal{L}(\lambda,\mu_v)]=[M(\mu_u):L(\mu_v)].\]

\emph{Step 3 (Kazhdan-Lusztig conjecture).} The weight $\mu_{\Id}+\rho=(a_1,\ldots,a_t)$ is dominant,
regular and integral, and $\mu_u=u^{-1}\cdot\mu_{\Id}$ for the dot action. By the Kazhdan-Lusztig
conjecture \cite{BB81, BK81}, in the normalization where $M(\mu_{\Id})$ is the dominant Verma module,
we have $[M(x\cdot\mu_{\Id}):L(y\cdot\mu_{\Id})]=P_{x,y}(1)$. Hence
$[M(\mu_u):L(\mu_v)]=P_{u^{-1},v^{-1}}(1)=P_{u,v}(1)$.
\end{proof}

By the inversion formula of Kazhdan and Lusztig (\cite{KL79}, (3.1)), it follows that for every $v\in S_t^0$
\begin{equation}\label{eq:inverseKL} [L(\caM_v)]=\sum_{u\in S_t^0,\ u\geq
v}(-1)^{l(u)-l(v)}P_{w_0u,w_0v}(1)\,[\caI_u].\end{equation}
For $v=\Id$ this is Tadi\'c's formula, since $P_{x,w_0}=1$ for all $x$.

\begin{lemma}\label{lem:ck} Let $1\leq k\leq t$.
\begin{enumerate}
\item[(i)] For $u\in S_t$, we have $u\geq c_k$ if and only if $u(1)\geq k$.
\item[(ii)] $P_{x,w_0c_k}=1$ for every $x\leq w_0c_k$.
\end{enumerate}
\end{lemma}

\begin{proof}
(i) We use the tableau criterion (\cite{BB05}, Theorem 2.6.3 and Proposition 2.4.8): $u\geq c_k$
if and only if, for every $j$, the increasingly ordered set $u(\{1,\ldots,j\})$ dominates
$c_k(\{1,\ldots,j\})$ entrywise.
 For $j\geq k$ we have $c_k(\{1,\ldots,j\})=\{1,\ldots,j\}$, and the condition is automatic.
  For $j<k$, $c_k(\{1,\ldots,j\})=\{1,\ldots,j-1,k\}$, and the condition reduces to $\max
  u(\{1,\ldots,j\})\geq k$.
   All these conditions follow from the case $j=1$, which is $u(1)\geq k$.

(ii) By (i), the interval $[e,w_0c_k]=w_0\,[c_k,w_0]$ is the set of $x\in S_t$ with $x(1)\leq t+1-k$.
This set is stable under right multiplication by the parabolic subgroup $W_J\simeq S_{t-1}$ fixing $1$.
Hence $w_0c_k$ is the longest element of its coset $w_0c_kW_J$, and the Schubert variety $X_{w_0c_k}$
is the inverse image, under the projection of the flag variety of $\mathbb{C}^t$ onto $\mathbb{P}^{t-1}$,
of a Schubert variety of $\mathbb{P}^{t-1}$. The latter is a linear subspace. Hence $X_{w_0c_k}$ is smooth,
and $P_{x,w_0c_k}=1$ for all $x\leq w_0c_k$ by \cite{KL80}.
\end{proof}

\begin{proof}[Second proof of Proposition \ref{conj:key}]
Let $A_k$ be as in Remark \ref{rmk:tadic_ck}. As in the proof of Proposition \ref{conj:key},
$A_k=\sum_{w\in S^0_t,\ w(1)\geq k}(-1)^{l(w)-l(c_k)}[\caI_w]$. If $\pi_k\neq 0$, then $c_k\in
S^0_t$, and by Lemma \ref{lem:ck} the condition $w(1)\geq k$ is equivalent to $w\geq c_k$, and
$P_{w_0w,w_0c_k}(1)=1$ for all such $w$. Comparing with (\ref{eq:inverseKL}) for $v=c_k$ gives
$A_k=[L(\caM_{c_k})]$. If $\pi_k=0$, then $A_k=0$ by Lemma \ref{lem:minors} (ii). Since
$[\pi_k\times\Delta'_k]=A_k+A_{k+1}$, the proposition follows.
\end{proof}

Similarly, by Proposition \ref{fact:KL}, the multiplicity in Proposition \ref{prop:multone} is
$P_{x,y}(1)$, and $P_{x,y}=1$ since $\deg P_{x,y}\leq(l(y)-l(x)-1)/2<1$.

\subsection{A geometric interpretation: Cousin resolutions}\label{ssec:geomint}
The proof of Theorem \ref{BGG} given above is combinatorial. We explain now why such a coresolution
exists, in the spirit of the geometric construction of the dual BGG resolution recalled in the
introduction. The point is the following criterion, which is a $p$-adic counterpart of the Cousin
resolution of the constant sheaf of a smooth variety.

\begin{prop}\label{prop:cousin} Let $v\in S_t^0$, and assume that $P_{w_0u,w_0v}=1$ for every $u\in
S_t^0$ with $u\geq v$. Then there is an exact sequence
\[0\rightarrow L(\caM_v)\rightarrow\caJ_v\rightarrow\bigoplus_{u\in S_t^0,\ u\geq v,\
l(u)=l(v)+1}\caJ_u\rightarrow\cdots\rightarrow\bigoplus_{u\in S_t^0,\ u\geq v,\ l(u)=l(v)+k}\caJ_u
\rightarrow\cdots.\]
\end{prop}

\noindent For $v=\Id$, the hypothesis holds because $P_{x,w_0}=1$ for all $x$ (the Schubert
variety $X_{w_0}$ is the whole flag variety, which is smooth), and we recover Theorem \ref{BGG}. For
$v=c_k$, it holds by Lemma \ref{lem:ck} (ii) ($X_{w_0c_k}$ is smooth), and we get a coresolution of
$L(\caM_{c_k})$ by the $\caJ_u$, $u(1)\geq k$, whose Euler characteristic is the formula of Remark
\ref{rmk:tadic_ck}. In general the hypothesis says that the Schubert variety $X_{w_0v}$ is
rationally smooth along the cells $X_{w_0u}$, $u\in S_t^0$, $u\geq v$, and the coresolution is the
image of the Cousin complex of its constant sheaf. Here the Schubert varieties only enter through their Kazhdan-Lusztig polynomials: the
geometry is that of the orbits in Zelevinsky's space (Step 1 below), whose Kazhdan-Lusztig
polynomials along the ladder strata are those of $S_t$. In particular, what matters for ladders is
the smoothness of the full flag variety of $\mathrm{GL}_t$, and not that of the Schubert variety
$X_{w'_0}$ attached to the interval $S_t^0$; this is consistent with the construction of \cite{BC15},
where the coresolution is the image of the dual BGG resolution of a finite-dimensional
$\mathfrak{gl}_t$-module under the Arakawa--Suzuki functor.

\begin{example} In the running example, Proposition \ref{prop:cousin} applies to $v=\Id$, which gives
the coresolution of Example \ref{ex:run2}, and to $v=c_k$, $2\leq k\leq6$, which gives a
coresolution of $L(\caM_{c_k})$ by the $\caJ_u$, $u\in S_6^0$, $u(1)\geq k$, i.e.\ by the rows
$k,\ldots,6$ of the table of Example \ref{ex:run5}.
\end{example}

\begin{proof}[Outline of proof]
\emph{Step 1 (the geometric category).} As in the proof of Proposition \ref{fact:KL}, the Bernstein
block is equivalent to the category of modules over the graded Hecke algebra $\mathbb{H}_\ell$. By
\cite{CG97}, Chapter 8, the finite-dimensional $\mathbb{H}_\ell$-modules with the central character
of $\pi$ (generalized) are the finite-dimensional modules, on which $H^{>0}_{G_V}(\mathrm{pt})$ acts
nilpotently, over the extension algebra $A=\mathrm{Ext}^\bullet_{G_V}(\mathcal{L},\mathcal{L})$ of
the $G_V$-equivariant derived category of Zelevinsky's space $E_V$ of graded nilpotent endomorphisms
(\cite{Zel85}), where $\mathcal{L}$ is the sum of the simple perverse sheaves $\mathrm{IC}_u$ attached
to the $G_V$-orbits. The orbits are the $O_u$ attached to the multisegments $\caM_u$, the stabilizers
are connected (they are open in the endomorphism algebras of representations of a quiver of type
$A$), $O_u\subset\overline{O_{u'}}$ if and only if $u\leq u'$, and $L(\caM_u)$ corresponds to the
simple module $L_u$. The module $\caJ_u$ corresponds to the graded dual $K_u^*$ of Lusztig's standard
module $K_u$, and the graded multiplicities are the Kazhdan--Lusztig polynomials,
$[K_u:L_{u'}]_q=P_{u,u'}(q)$ (the graded form of Proposition \ref{fact:KL}).

\emph{Step 2 (Kato's structure theorem).} The conditions of \cite{Kat17}, \S1, hold in this
situation: the stabilizers are connected and the $\mathrm{IC}_u$ are pointwise pure (Lusztig).
For finitely generated graded $A$-modules $M,N$, we write, following \cite{Kat17}, \S1,
$\mathrm{Ext}^i_A(M,N)=\bigoplus_j\mathrm{Ext}^i(M,N\langle j\rangle)$, the extension groups being computed in
the category of graded $A$-modules and $N\langle j\rangle$ denoting the shift of the grading; it is a graded
vector space, and for a graded vector space $E=\bigoplus_jE_j$ with finite-dimensional pieces we write
$\mathrm{gdim}\,E=\sum_j\dim(E_j)\,q^{j/2}$ for its graded dimension. Hence (\cite{Kat17}, Theorem C) the modules $\widetilde K_u=P_u/\sum_{u'<u}Ae_{u'}P_u$ satisfy
$\mathrm{Ext}^i_A(\widetilde K_u,K_{u'}^*)=\mathbb{C}$ if $u=u'$ and $i=0$, and $0$ otherwise, and
$[P_u:\widetilde K_{u'}]=[K_{u'}:L_u]$.

\emph{Step 3 (the Ext groups).} By Step 2, in the graded Grothendieck group,
\[\textstyle\sum_i(-1)^i\,\mathrm{gdim}\,\mathrm{Ext}^i(\widetilde K_u,L_v)\]
is the coefficient of $[K_u^*]$ in $[L_v]$, that is,
the $(u,v)$ entry of the inverse of the matrix $(P_{u,u'}(q))$. By the inversion formula of Kazhdan and
Lusztig (\cite{KL79}, (3.1)), this entry is $(-1)^{l(u)-l(v)}P_{w_0u,w_0v}$, up to a power of $q$. By
purity, $\mathrm{Ext}^i(\widetilde K_u,L_v)$ is pure of weight $i$, so that there is no cancellation,
and the hypothesis gives $\mathrm{Ext}^i(\widetilde K_u,L_v)=\mathbb{C}$ if $u\geq v$ and
$i=l(u)-l(v)$, and $0$ otherwise.

\emph{Step 4 (the coresolution).} Put $Y_0=L_v$, and inductively let $f_k:Y_k\to\bigoplus_{l(u)=l(v)+k}
K_u^*$ be the map whose components span the one-dimensional spaces $\mathrm{Hom}(Y_k,K_u^*)$, and
$Y_{k+1}=\mathrm{coker}(f_k)$. By Step 2, $\mathrm{Ext}^i(\widetilde K_u,Y_k)$ is obtained from
$\mathrm{Ext}^{i+k}(\widetilde K_u,L_v)$, so that it vanishes unless $i=l(u)-l(v)-k$; this gives the
injectivity of $f_k$ and the exactness, as in Johnson's construction \cite{Joh84}. We omit the details, which are standard: the dictionary of
Step 1, including the graded form of Proposition \ref{fact:KL}, follows from \cite{CG97}, \S8.6; the
purity in Step 3 follows from the parity vanishing and the purity of the stalks and costalks of the
$\mathrm{IC}_u$; and Step 4 is the usual construction of linear coresolutions in a graded highest
weight category.
\end{proof}

\begin{lemma}\label{lem:dimorb} Let $V=\bigoplus_xV_x$ be the graded vector space attached to the
cuspidal support of $\pi$, $G_V=\prod_x\mathrm{GL}(V_x)$, and $O_u$ the $G_V$-orbit in Zelevinsky's
space $E_V$ attached to the multisegment $\caM_u$ (\cite{Zel85}, Proposition 1). For $u\in S_t^0$,
\[\dim O_u=\dim O_{\Id}+l(u).\]
\end{lemma}

\begin{proof}
Let $N\in O_u$. Its Jordan strings are the non-empty segments $s_i=[a_{u(i)};b_i]$, $s_i$ running from
degree $b_i$ down to degree $a_{u(i)}$. We have $\dim O_u=\dim G_V-\dim Z(N)$, where $\dim G_V$ does not
depend on $u$, and $Z(N)$ is the space of degree-preserving endomorphisms of $V$ commuting with $N$.
A degree-preserving map from the string $s$ to the string $s'$ commuting with $N$ is determined by the
image of the top vector of $s$, which must be killed by the power of $N$ of exponent the length of $s$;
so these maps form a space of dimension $1$ if the bottom and top degrees satisfy
$a(s)\leq a(s')\leq b(s)\leq b(s')$, and $0$ otherwise. Since $b_1>\cdots>b_t$ and $a_1>\cdots>a_t$, for
$i<j$ only the maps from $s_j$ to $s_i$ can occur, and they do if and only if $u(i)<u(j)$ and
$a_{u(i)}\leq b_j$ (both segments are then non-empty, since $a_{u(j)}<a_{u(i)}\leq b_j<b_i$). Hence
\[\dim Z(N)=n(u)+\#\{i<j:u(i)<u(j),\ a_{u(i)}\leq b_j\},\]
where $n(u)$ is the number of non-empty segments. If $i<j$ and $u(i)>u(j)$, then $a_{u(i)}<a_{u(j)}\leq
b_j+1$, so $a_{u(i)}\leq b_j$. So the second term is $\#\{i<j:a_{u(i)}\leq b_j\}-l(u)$. Let
$c(k)=\#\{j:b_j\geq a_k\}$. Then $\#\{j>i:a_{u(i)}\leq b_j\}=\max(c(u(i))-i,0)$, and $s_i$ is non-empty if and
only if $i\leq c(u(i))$. Moreover $c(u(i))\geq i-1$: this is clear for $i=1$, and for $i>1$,
$b_{i-1}>b_i\geq a_{u(i)}-1$. Hence
\[n(u)+\#\{i<j:a_{u(i)}\leq b_j\}=\sum_{i=1}^t\big(c(u(i))-i+1\big)=\sum_{k=1}^tc(k)-\frac{t(t-1)}2,\]
which does not depend on $u$, and $\dim Z(N)=\dim Z(N_{\Id})-l(u)$.
\end{proof}

\subsection{Balanced representations}\label{ssec:balanced}
The geometric argument of Proposition \ref{prop:cousin} does not use the specific shape of ladders,
and Theorem \ref{BGG} extends to the class of \emph{balanced} representations introduced by
Lapid and M\'inguez \cite{LM18}. We only sketch this extension, leaving the details to the reader.

Let $\frmm=[a_1;b_1]+\cdots+[a_k;b_k]$ be a multisegment with $b_1>\cdots>b_k$ and $a_1,\ldots,a_k$
pairwise distinct; following \cite{LM18}, such a multisegment is called \emph{regular}. Let
$a'_1<\cdots<a'_k$ be the beginnings in increasing order and, for $\tau\in S_k$, put
\[\frmm_\tau=[a'_{\tau^{-1}(1)};b_1]+\cdots+[a'_{\tau^{-1}(k)};b_k],\]
with the conventions of the introduction for the segments $[a;a-1]$ (empty) and $[a;b]$, $b<a-1$
(the corresponding representations are $0$). Then $\frmm=\frmm_\sigma$ for a unique $\sigma\in S_k$, and
there is a unique $\sigma_0\in S_k$ such that $\frmm_\tau$ has no segment $[a;b]$ with $b<a-1$
(i.e.\ the attached representations are non-zero) if and only if $\tau\geq\sigma_0$ (\cite{LM18},
\S6.1); moreover
$\sigma_0\leq\sigma$, and for $\tau,\tau'\geq\sigma_0$, $L(\frmm_{\tau'})$ is a subquotient of the
standard module attached to $\frmm_\tau$ if and only if $\tau\leq\tau'$. For a ladder, $\sigma$ is the
longest element of $S_k$, and the multisegments $\frmm_\tau$, $\tau\in[\sigma_0,\sigma]$, are the
$\caM_w$, $w\in S_t^0$.

Lapid and M\'inguez call $\frmm$ \emph{balanced} when a combinatorial condition on its degenerations
holds (two numbers attached to $\frmm$, its depth and its complexity, coincide; \cite{LM18},
Definition 6.5). By \cite{LM18}, Theorems 1.2 and 7.1, for a regular $\frmm$ the following conditions
are equivalent: $\frmm$ is balanced; $\langle\frmm\rangle\times\langle\frmm\rangle$ is irreducible;
$P_{\tau,\sigma}=1$ for every $\tau\in[\sigma_0,\sigma]$, i.e.\ the Schubert variety $X_\sigma$ is
rationally smooth along the cells attached to $[\sigma_0,\sigma]$; and Tadi\'c's formula holds in the
form
\[\langle\frmm\rangle=\sum_{\tau\in[\sigma_0,\sigma]}\mathrm{sgn}(\sigma\tau)\,
\langle[a'_{\tau^{-1}(1)};b_1]\rangle\times\cdots\times\langle[a'_{\tau^{-1}(k)};b_k]\rangle\]
in the Grothendieck group. Applying the Zelevinsky involution, which maps $\langle\frn\rangle$ to
$L(\frn)$ and $\langle\Delta\rangle$ to $\Delta(\Delta)$, the last condition expresses $L(\frmm)$ as the
alternating sum of the costandard modules $\caJ(\frmm_\tau)$ attached to the $\frmm_\tau$,
$\tau\in[\sigma_0,\sigma]$ (products of the $\Delta([a'_{\tau^{-1}(i)};b_i])$ written in increasing
order of their ends), with the sign $(-1)^{l(\sigma)-l(\tau)}$. For ladders, this is Theorem
\ref{LM}. Moreover, the graded multiplicities of the irreducible representations $L(\frmm_{\tau'})$ in
the $\caJ(\frmm_\tau)$ are the Kazhdan--Lusztig polynomials $P_{\tau,\tau'}$ of $S_k$ (\cite{LM18},
Corollary 10.1, which uses the functor of Arakawa--Suzuki; geometrically, the closures of the orbits
attached to the $\frmm_\tau$ in Zelevinsky's space are locally isomorphic to Schubert varieties,
\cite{Zel85}).

\begin{rmk}\label{rmk:balanced} Let $\frmm$ be a regular balanced multisegment. Then $L(\frmm)$
has a coresolution
\[0\rightarrow L(\frmm)\rightarrow\caJ(\frmm_\sigma)\rightarrow\bigoplus_{\substack{\tau\in[\sigma_0,\sigma]\\
l(\tau)=l(\sigma)-1}}\caJ(\frmm_\tau)\rightarrow\cdots\rightarrow\bigoplus_{\substack{\tau\in[\sigma_0,\sigma]\\
l(\tau)=l(\sigma)-k}}\caJ(\frmm_\tau)\rightarrow\cdots\rightarrow\caJ(\frmm_{\sigma_0})\rightarrow0,\]
whose Euler characteristic is the formula above. The argument is that of Proposition
\ref{prop:cousin}. Steps 1 and 2 only use the description of the block as modules over the extension
algebra of Zelevinsky's space, with connected stabilizers. In Step 3, the Kazhdan--Lusztig
polynomials of $S_t$ are replaced by those of $S_k$ attached to $[\sigma_0,\sigma]$, and the
hypothesis $P_{\tau,\sigma}=1$ is exactly the balanced condition: as for ladders, the coresolution is
the image of the Cousin complex of the constant sheaf of the closure of the orbit of $\frmm$, which is
rationally smooth along the orbits of the $\frmm_\tau$, and the degrees are the codimensions of these
orbits, which are the differences of lengths. Step 4 is unchanged. The maps of the coresolution are
determined up to scalars, since the relevant spaces of homomorphisms are one-dimensional, but we do not
have an explicit construction as in Section \ref{sec:complex}, which relies on the minors of ladders.
\end{rmk}

% Bibliography added on 23/10/2026: volume and page numbers to be checked.

\end{document}